\documentclass[reqno,12pt]{amsart} 
\usepackage{amsmath}
\usepackage{amssymb}
\usepackage{amsthm}
\newcommand\NoBlackBoxes{\global\overfullrule0pt}
\NoBlackBoxes
\theoremstyle{plain} 

\def\4{\kern1pt}

\def\6{\vphantom0}

\def\8{\kern-10pt}
\def\7#1{_{(#1)}}

\makeatletter
\let\serieslogo@\relax
\let\@setcopyright\relax

\def\speciallabelmark#1{\def\@currentlabel{#1}}
\makeatother

\begin{document}

\def\ffrac#1#2{\raise.5pt\hbox{\small$\4\displaystyle\frac{\,#1\,}{\,#2\,}\4$}}
\def\ovln#1{\,{\overline{\!#1}}}
\def\ve{\varepsilon}
\def\kar{\beta_r}

\title{FISHER INFORMATION AND THE CENTRAL LIMIT THEOREM
}

\author{S. G. Bobkov$^{1,4}$}
\thanks{1) School of Mathematics, University of Minnesota, USA;
Email: bobkov@math.umn.edu}
\address
{Sergey G. Bobkov \newline
School of Mathematics, University of Minnesota  \newline 
127 Vincent Hall, 206 Church St. S.E., Minneapolis, MN 55455 USA
\smallskip}
\email {bobkov@math.umn.edu} 

\author{G. P. Chistyakov$^{2,4}$}
\thanks{2) Faculty of Mathematics, University of Bielefeld, Germany;
Email: chistyak@math.uni-bielefeld.de}
\address
{Gennadiy P. Chistyakov\newline
Fakult\"at f\"ur Mathematik, Universit\"at Bielefeld\newline
Postfach 100131, 33501 Bielefeld, Germany}
\email {chistyak@math.uni-bielefeld.de}

\author{F. G\"otze$^{3,4}$}
\thanks{3) Faculty of Mathematics, University of Bielefeld, Germany;
Email: goetze@math.uni-bielefeld.de}
\thanks{4) Research partially supported by 
NSF grant DMS-1106530 and SFB 701}
\address
{Friedrich G\"otze\newline
Fakult\"at f\"ur Mathematik, Universit\"at Bielefeld\newline
Postfach 100131, 33501 Bielefeld, Germany}
\email {goetze@mathematik.uni-bielefeld.de}



\begin{abstract}
An Edgeworth-type expansion is established for the relative Fisher information 
distance to the class of normal distributions of sums of i.i.d. random 
variables, satisfying moment conditions. The validity of the central limit 
theorem is studied via properties of the Fisher information along convolutions.
\end{abstract}

\maketitle
\markboth{S. G. Bobkov, G. P. Chistyakov and F. G\"otze}{Fisher Information}




\def\theequation{\thesection.\arabic{equation}}
\def\E{{\bf E}}
\def\R{{\bf R}}
\def\C{{\bf C}}
\def\P{{\bf P}}
\def\H{{\rm H}}
\def\Im{{\rm Im}}
\def\Tr{{\rm Tr}}

\def\k{{\kappa}}
\def\M{{\cal M}}
\def\Var{{\rm Var}}
\def\Ent{{\rm Ent}}
\def\O{{\rm Osc}_\mu}

\def\ep{\varepsilon}
\def\phi{\varphi}
\def\F{{\cal F}}
\def\L{{\cal L}}

\def\be{\begin{equation}}
\def\en{\end{equation}}
\def\bee{\begin{eqnarray*}}
\def\ene{\end{eqnarray*}}


\section{{\bf Introduction}}
\setcounter{equation}{0}

\vskip2mm
\noindent
Given a random variable $X$ with an absolutely continuous density $p$,
the Fisher information of $X$ (or its distribution) is defined by
$$
I(X) = I(p) = \int_{-\infty}^{+\infty} \frac{p'(x)^2} {p(x)} \, dx,
$$
where $p'$ denotes a Radon-Nikodym derivative of $p$. In all other cases,
let $I(X) = +\infty$.

With the first two moments of $X$ being fixed, this quantity is minimized for 
the normal distribution (which is a variant of Cram\'er-Rao's inequality). 
That is, if $\E X = a$, $\Var(X) = \sigma^2$, then we have $I(X) \geq I(Z)$ 
for $Z \sim N(a,\sigma^2)$ with density
$$
\varphi_{a,\sigma}(x) = 
\frac{1}{\sqrt{2 \pi \sigma^2}}\ e^{-(x-a)^2/2\sigma^2}. 
$$
Moreover, the equality $I(X) = I(Z)$ holds if and only if $X$ is normal. 

In many applications the relative Fisher information
$$
I(X||Z) = I(X)- I(Z) = \int_{-\infty}^{+\infty} \bigg(\frac{p'(x)}{p(x)} -
\frac{\varphi_{a,\sigma}'(x)}{\varphi_{a,\sigma}(x)}\bigg)^2\, p(x)\, dx,
$$
which is used as a strong measure of non-Gaussianity of $X$. For example, it 
dominates the relative entropy, or Kullback-Leibler distance of the 
distribution of $X$ to the standard normal distribution; more precisely 
(cf. Stam [S]),
\be
\frac{\sigma^2}{2}\,I(X||Z) \, \geq \, D(X||Z) = \,\int_{-\infty}^{+\infty} 
p(x) \log \frac{p(x)}{\varphi_{a,\sigma}(x)}\ dx.
\en

We consider the scheme of a sequence of sums of independent identically 
distributed random variables $(X_n)_{n \geq 1}$. Assuming that $\E X_1 = 0$, 
$\Var(X_1) = 1$, define the normalized sums
$$
Z_n = \frac{X_1 + \dots + X_n}{\sqrt{n}}.
$$
Since $Z_n$ are weakly convergent in distribution to $Z \sim N(0,1)$, one may 
wonder whether the convergence holds in a stronger sense. A remarkable 
observation in this respect is due to Barron and Johnson proving in [B-J] 
that 
\be
I(Z_n) \rightarrow I(Z), \quad {\rm as} \ \ n \rightarrow \infty,
\en
i.e., $I(Z_n||Z) \rightarrow 0$, if and only if $I(Z_{n_0})$ is finite 
for some $n_0$. In particular, it suffices to require that $I(X_1) < +\infty$, 
although choosing larger values of $n_0$ considerably enhances 
the range of applicability of this theorem.

Quantitative estimates on the relative Fisher information in the central limit 
theorem are partly developed, as well. In the i.i.d. case Barron and 
Johnson [B-J], and Artstein, Ball, Barthe and Naor [A-B-B-N1] derived 
an asymptotic bound $I(Z_n||Z) = O(1/n)$ under the hypothesis that the 
distribution of $X_1$ admits an analytic inequality of Poincar\'e-type 
(cf. also [J]). Poincar\'e inequalities involve a large variety of "nice" 
probability distributions on the line all having finite exponential moments.

One of the aims of this paper is to study the exact asymptotics (or rates) of 
$I(Z_n||Z)$ under standard moment conditions. We prove:

\vskip5mm
{\bf Theorem 1.1.} {\it Let $\E\, |X_1|^s < +\infty$ for an integer $s \geq 2$, 
and assume $I(Z_{n_0}) < +\infty$, for some $n_0$. Then for certain 
coefficients $c_j$ we have, as $n \rightarrow \infty$,
\be
I(Z_n||Z) = \frac{c_1}{n} + \frac{c_2}{n^2} + \dots + 
\frac{c_{[(s-2)/2)]}}{n^{[(s-2)/2)]}} + 
o\left(n^{-\frac{s-2}{2}} \, (\log n)^{-\frac{(s-3)_+}{2}}\right).
\en
}

\vskip2mm
As it turns out, a similar expansion holds as well for the entropic distance
$D(Z_n||Z)$, cf. [B-C-G2], showing a number of interesting analogies in the 
asymptotic behavior of these two distances. In particular, in both cases 
each coefficient $c_j$ is given by a certain polynomial in the cumulants
$\gamma_3,\dots,\gamma_{2j+1}$. 

In order to describe these polynomials, we first 
note that, by the moment assumption, the cumulants
$$
\gamma_r = i^{-r}\, \frac{d^r}{dt^r} \, \log \E\, e^{itX_1}|_{t=0}
$$ 
are well-defined for all positive integers $r \leq s$, and one may introduce 
the well-known functions
$$
q_k(x) \ = \, \varphi(x)\, \sum H_{k + 2j}(x) \,
\frac{1}{r_1!\dots r_k!}\, \bigg(\frac{\gamma_3}{3!}\bigg)^{r_1} \dots
\bigg(\frac{\gamma_{k+2}}{(k+2)!}\bigg)^{r_k}
$$
involving the Chebyshev-Hermite polynomials $H_k$. Here 
$\varphi = \varphi_{0,1}$ denotes the density of the standard normal law, 
and the summation runs over all non-negative integer solutions 
$(r_1,\dots,r_k)$ to the equation
$r_1 + 2 r_2 + \dots + k r_k = k$ with $j = r_1 + \dots + r_k$.

The functions $q_k$ are correctly defined for $k = 1,\dots,s-2$. They appear 
in Edgeworth-type expansions approximating the density of $Z_n$. We shall 
employ them to derive an expansion in powers of $1/n$ for the 
distance $I(Z_n||Z)$, which leads us to the following description
of the coefficients in (1.3),
\be
c_j \, = \, \sum_{k=2}^{2j}\, (-1)^k \sum \int_{-\infty}^{+\infty}
(q_{r_1}' + x q_{r_1}) (q_{r_2}' + x q_{r_2})\,
q_{r_3} \dots q_{r_k}\, \frac{dx}{\varphi^{k - 1}}.
\en
Here, the inner summation is carried out over all positive integer tuples
$(r_1,\dots,r_k)$ such that $r_1 + \dots + r_k = 2j$.

For example, $c_1 = \frac{1}{2}\,\gamma_3^2$, and in the case $s=4$
(1.3) becomes
\be
I(Z_n||Z) = \frac{1}{2 n}\,\left(\E X_1^3\right)^2 + 
o\left(\frac{1}{n\, (\log n)^{1/2}}\right).
\en
Hence, under the 4-th moment condition, we have $I(Z_n||Z) \leq \frac{C}{n}$
with some constant $C$ (which can actually be chosen to depend on 
$\E X_1^4$ and $I(X_1)$, only).

For $s=6$, the result involves the coefficient $c_2$ which depends on
$\gamma_3,\gamma_4$, and $\gamma_5$. If $\gamma_3 = 0$ (i.e. $\E X_1^3 = 0$), 
we have $c_1 = 0$, $c_2 = \frac{1}{6}\,\gamma_4^2$, and then
$$
I(Z_n||Z) = \frac{1}{6 n^2}\,\left(\E X_1^4 - 3\right)^2 + 
o\left(\frac{1}{n^2\, (\log n)^{3/2}}\right).
$$

More generally, the representation (1.3) simplifies, if the first $k-1$
moments of $X_1$ coincide with the corresponding moments of $Z \sim N(0,1)$.

\vskip5mm
{\bf Corollary 1.2.} {\it Let $\E\, |X_1|^s < +\infty$ $(s \geq 4)$, and assume
$I(Z_{n_0}) < +\infty$, for some $n_0$. Given $k = 3,4,\dots,s$, assume that
$\gamma_j = 0$ for all $3 \leq j < k$. Then
\be
I(Z_n||Z) \, = \, \frac{\gamma_k^2}{(k-1)!}\cdot \frac{1}{n^{k-2}} + 
O\bigg(\frac{1}{n^{k-1}}\bigg) +
o\bigg(\frac{1}{n^{(s-2)/2}\, (\log n)^{(s-3)/2}}\bigg).
\en
}

\vskip2mm
This relation is consistent with an observation of Johnson who noticed that 
if $\gamma_k \neq 0$, $I(Z_n||Z)$ cannot be 
asymptotically better than $n^{-(k-2)}$ ([J], Lemma 2.12).

Note that if $k < \frac{s}{2}$, the $O$-term in (1.6) dominates the $o$-term. 
But when $k \geq \frac{s}{2}$ it can be removed, and if $k > \frac{s}{2} + 1$, 
(1.6) just says that
\be
I(Z_n||Z) \, = \, o\big(n^{-(s-2)/2} \, (\log n)^{-(s-3)/2}\big).
\en

For the values $s=2,3$ there are no coefficients $c_j$ in the sum (1.3).
In case $s=2$ Theorem 1.1 reduces to Barron-Johnson's theorem (1.2), while
under a 3-rd moment assumption we only have
$$
I(Z_n||Z) = o\Big(\frac{1}{\sqrt{n}}\Big).
$$
A similar observation holds for the whole range of reals $2<s<4$. 
Here the expansion (1.3) should be replaced by the bound (1.7).
Although this bound is worse than (1.5), it cannot be essentially 
improved. As shown in [B-C-G2], it may happen that $\E\, |X_1|^s < +\infty$ 
with $D(X_1) < +\infty$ (in fact, with $I(X_1) < +\infty$), while
$$
D(Z_n||Z)\, \geq\, \frac{c}{n^{(s-2)/2} \ (\log n)^\eta},
\qquad n \geq n_1(X_1),
$$
where the constant $c>0$ depends on $s$ and an arbitrary prescribed value 
$\eta > s/2$. In view of (1.1), a similar lower bound therefore 
holds for $I(Z_n||Z)$, as well.

Another interesting issue connected with the convergence theorem (1.2) and 
the expansion (1.3) is the characterization of distributions for which these 
results hold. Indeed, the condition $I(X_1) < +\infty$ corresponding to 
$n_0 = 1$ in Theorem 1.1 seems to be way too strong. To this aim, 
we establish an explicit criterion such that $I(Z_{n_0}) < +\infty$ holds
for sufficiently large $n_0$ in terms of the characteristic function
$f_1(t) = \E\,e^{itX_1}$ of $X_1$.

\vskip5mm
{\bf Theorem 1.3.} {\it Given independent identically distributed 
random variables $(X_n)_{n \geq 1}$ with finite second moment, 
the following assertions are equivalent{\rm :}

\vskip2mm
$a)$ \ For some $n_0$, $Z_{n_0}$ has finite Fisher information{\rm ;}

$b)$ \ For some $n_0$, $Z_{n_0}$ has density of bounded total variation{\rm ;}

$c)$ \ For some $n_0$, $Z_{n_0}$ has a continuously differentiable density 
$p_{n_0}$ such that
$$
\int_{-\infty}^{+\infty} |p_{n_0}'(x)|\, dx < +\infty;
$$

$d)$ \ For some $\ep>0$, $|f_1(t)| = O(t^{-\ep})$, as $t \rightarrow +\infty${\rm ;}

$e)$ \ For some $\nu > 0$,
\be
\int_{-\infty}^{+\infty} |f_1(t)|^\nu\,|t|\, dt < +\infty.
\en
}

Property $c)$ is a formally strengthened variant of $b)$, although in general
they are not equivalent. (For example, the uniform distribution has density 
of bounded total variation, but its density is not everywhere differentiable.)

Properties $a)-c)$ are equivalent to each other without any moment assumption,
while $d)-e)$ are always necessary for the finiteness of $I(Z_n)$ with 
large $n$. These two last conditions show that the range of applicability 
of Theorem 1.1 is indeed rather wide, since almost all reasonable 
absolutely continuous distributions satisfy (1.8). 
The latter should be compared to and viewed as a certain strengthening 
of the following condition (sometimes called a smoothness condition)
$$
\int_{-\infty}^{+\infty} |f_1(t)|^\nu\, dt < +\infty, \qquad {\rm for \ some} \ \
\nu > 0.
$$
It is equivalent to the property that, for some $n_0$, $Z_{n_0}$ has 
a bounded continuous density $p_{n_0}$ (cf. e.g. [BR-R]). In this and only
in this case, a uniform local limit theorem holds: 
$\Delta_n = \sup_x |p_n(x) - \varphi(x)| \rightarrow 0$, as $n \rightarrow \infty$.
That this assertion is weaker compared to the convergence in Fisher 
information distance such as (1.2) can be seen by Shimizu's inequality 
$\Delta_n^2 \leq c I(Z_n||Z)$, which holds with some absolute constant $c$
([Sh], [B-J], Lemma 1.5). Note in this connection that Shimizu's 
inequality may be strengthened in terms of the total variation distance as
$\|p_n - \varphi\|_{\rm TV}^2 \leq c I(Z_n||Z)$. Using Theorem 1.3, 
this shows that (1.2) is equivalent to the convergence
$\|p_n - \varphi\|_{\rm TV} \rightarrow 0$.

The paper is organized in the following way.
We start with the description of general properties of densities having
finite Fisher information (Section 2) and properties of Fisher information
as a functional on spaces of densities (showing lower semi-continuity and 
convexity, Section 3). Some of the properties and relations which we state 
for completeness may be known already. We apologize for being unable to find 
references for them.

In Sections 4-5 we turn to upper bounds needed mainly in the proof of 
Theorem 1.3. Further properties of densities emerging after several
convolutions, as well as, bounds under additional moment assumptions are 
discussed in Sections 6-8. In Section 9 we complete the proof of Theorem 1.3, 
and in the next section we state basic lemmas on Edgeworth-type expansions 
which are needed in the proof of Theorem 1.1. 
Sections 11-12 are devoted to the proof itself.
Some remarks leading to the particular case $s=2$ in Theorem 1.1
(Barron-Johnson theorem)
are given in Section 13. Finally, in the last section we briefly describe
the modifications needed to obtain Theorem 1.1 under moment assumptions
with arbitrary real values of $s$.

\vskip1mm
{\it Table of contents}

{\small
1. Introduction

2. General properties of densities with finite Fisher information

3. Fisher information as a functional

4. Convolution of three densities of bounded variation

5. Bounds in terms of characteristic functions

6. Classes of densities representable as convolutions

7. Bounds under moment assumptions

8. Fisher information in terms of the second derivative

9. Normalized sums. Proof of Theorem 1.3

10. Edgeworth-type expansions

11. Behavior of densities not far from the origin

12. Moderate deviations

13. Theorem 1.1 in the case $s=2$ and Corollary 1.2

14. Extensions to non-integer $s$. Remarks on lower bounds

}


\vskip10mm
\section{{\bf General properties of densities with finite Fisher information}}
\setcounter{equation}{0}

\vskip2mm
\noindent
If a random variable $X$ has density $p$ with finite Fisher information
\be
I(X) = I(p) = \int_{-\infty}^{+\infty} \frac{p'(x)^2} {p(x)} \, dx,
\en
$p$ has to be absolutely continuous, and then the derivative $p'(x)$ 
exists and is finite on a set of full Lebesgue measure. 

One may write an equivalent definition by involving the score function
$\rho(x) = \frac{p'(x)}{p(x)}$. In general $\P\{p(X) > 0\} = 1$, so
the random variable $\rho(X)$ is well defined with probability~1, and thus
\be
I(X) = \E\,\rho(X)^2.
\en
However, strictly speaking, the integration in (2.1) should be restricted to 
the open set $\{x: p(x) > 0\}$.

For different purposes, it is useful to realize how the ratio 
$\frac{p'(x)^2}{p(x)}$ may behave when $p(x)$ is small and is even vanishing. 
The behavior cannot be arbitrary, when the Fisher information is finite.
The following statement plays a "justifying" role in obtaining of many Fisher 
information bounds on the density and its derivatives.

\vskip5mm
{\bf Proposition 2.1.} {\it Assume $X$ has density $p$ with finite Fisher 
information. If $p$ is differentiable at the point $x_0$ such that
$p(x_0) = 0$, then $p'(x_0) = 0$.
}

\vskip5mm
{\bf Proof.}
If $p$ is differentiable in some neighborhood of $x_0$ and its
derivative is continuous at this point, the statement is obvious.

To cover the general case, for simplicity of notations let $x_0 = 0$ and assume 
that $c = p'(0) > 0$. Since $p(\ep) = c\ep + o(\ep)$, as $\ep \rightarrow 0$,
one may choose $\ep_0>0$ 
such that
$$
\frac{3c}{4}\,|x| \leq p(x) \leq \frac{5c}{4}\,|x|, \quad {\rm for \ all} \ 
0 \leq |x| \leq \ep_0.
$$
In particular, $p$ is positive on $(0,\ep_0]$. Hence, by the definition (2.1),
$$
I(X) \geq \int_0^{\ep_0} \frac{p'(x)^2}{p(x)}\,dx \geq
\frac{4}{5c}\, \int_0^{\ep_0} \frac{p'(x)^2}{x}\,dx.
$$
We split the last integral into the intervals 
$\Delta_n = (2^{-(n+1)}\ep_0,2^{-n}\ep_0)$ and then estimate $p(x)$ from above
on each of them, which leads to
$$
\frac{5c \ep_0}{4}\,I(X) \, \geq \, \sum_{n=0}^\infty \, 
2^n \int_{\Delta_n} p'(x)^2\,dx.
$$
Now, applying Cauchy's inequality and using 
$p(x) - p(\frac{x}{2}) \geq \frac{c}{8}\,x$ for $0 \leq x \leq \ep_0$,
we obtain
\bee
\int_{\Delta_n} p'(x)^2\,dx 
 & \geq &
2^{n+1} \bigg(\int_{\Delta_n} p'(x)\,dx\bigg)^2 \\
 & = &
2^{n+1}\, \left(p\big(2^{-n}\ep_0\big) - p\big(2^{-(n+1)}\ep_0\big)\right)^2
\, \geq \, 2^{-(n+1)}\,\frac{(c \ep_0)^2}{64}.
\ene
As a result,
$$
\frac{5c \ep_0}{4}\,I(X) \, \geq \, \sum_{n=0}^\infty 2^n \cdot 
2^{-(n+1)} \cdot \frac{(c \ep_0)^2}{64} = +\infty,
$$
a contradiction with finiteness of the Fisher information.
Proposition 2.1 is proved.

\vskip2mm
As an example illustrating a possible behavior as in Proposition 2.1, 
one may consider the beta distribution with parameters 
$\alpha = \beta = 3$, which has density
$$
p(x) = 30\,(x(1-x))^2, \qquad 0 \leq x \leq 1.
$$
Then $X$ has finite Fisher information, although $p(x_0) = p'(x_0) = 0$
at $x_0 = 0$ and $x_0 = 1$. 

More generally, if a density $p$ is supported and twice differentiable 
on a finite interval $[a,b]$, and if $p$ has finitely many zeros 
$x_0 \in [a,b]$, and $p'(x_0) = 0$, $p''(x_0) > 0$ 
at any such point, then $X$ has finite Fisher information.

Now, let us return to the definitions (2.1)-(2.2). By Cauchy's inequality,
$$
I(X)^{1/2} = \big(\E\,\rho(X)^2\big)^{1/2} \geq \E\,|\rho(X)| =
\int_{\{p(x)>0\}} |p'(x)|\,dx.
$$
Here, by Proposition 2.1, the last integral may be extended to the whole 
real line without any change, and then it represents the total variation
of the function $p$ in the usual sense of the Theory of Functions:
$$
\|p\|_{{\rm TV}} \, = \, \sup\, \sum_{k=1}^n |p(x_k) - p(x_{k-1})|,
$$
where the supremum runs over all finite collections 
$x_0 < x_1 < \dots < x_n$. 

In the sequel, we consider this norm also for densities which are not 
necessarily continuous, and then it is natural to require that, for each $x$, 
the value $p(x)$ lies in the closed segment $\Delta(x)$ with endpoints $p(x-)$ 
and $p(x+)$. Note that if we change $p(x)$ at a point of discontinuity 
such that $p(x)$ goes out of $\Delta(x)$, then the measure with density 
$p$ is unchanged, while $\|p\|_{\rm TV}$ will increase.

Thus, if the Fisher information $I(X)$ is finite, the density
$p$ of $X$ is a function of bounded variation, so the limits
$$
p(-\infty) = \lim_{x \rightarrow -\infty} p(x), \qquad
p(+\infty) = \lim_{x \rightarrow +\infty} p(x)
$$
exist and are finite. But, since $p$ is a density (hence integrable), these 
limits must be zero. In addition, for any $x$,
$$
p(x) = \int_{-\infty}^x p'(y)\,dy \leq \int_{-\infty}^x |p'(y)|\,dy \leq
\sqrt{I(X)}.
$$
We can summarize these elementary observations in the following:

\vskip5mm
{\bf Proposition 2.2.} {\it If $X$ has density $p$ with finite Fisher 
information $I(X)$, then $p(-\infty) = p(+\infty) = 0$, and the density 
has finite total variation satisfying
$$
\|p\|_{{\rm TV}} = \int_{-\infty}^{+\infty} |p'(x)|\,dx \leq \sqrt{I(X)}.
$$
In particular, $p$ is bounded: $\max_x p(x) \leq \sqrt{I(X)}$.
}

\vskip5mm
{\bf Corollary 2.3.} {\it If $X$ has finite Fisher information, then its 
characteristic function $f(t) = \E\,e^{itX}$ admits the bound
$$
|f(t)| \leq \frac{1}{|t|}\,\sqrt{I(X)}, \qquad t \in \R.
$$
}

\vskip2mm
Indeed, using Proposition 2.2, one may integrate by parts,
$$
it\, \E\,e^{itX} = \int_{-\infty}^{+\infty} p(x)\,d\,e^{itx} = 
-\int_{-\infty}^{+\infty} e^{itx}\,p'(x)\,dx,
$$
which gives
$|t|\, |\E\,e^{itX}| \leq \int_{-\infty}^{+\infty} |p'(x)|\,dx \leq \sqrt{I(X)}$.

\vskip2mm
Another immediate consequence of Proposition 2.2 is that both $p$ and $p'$ are 
square integrable, that is, they belong to the Sobolev space 
$W_1^2 = W_1^2(-\infty,+\infty)$ of all absolutely continuous functions on 
the real line with finite Euclidean (Hilbert) norm
$$
\|u\|_{W_1^2}^2 = \int_{-\infty}^{+\infty} u(x)^2\,dx +
\int_{-\infty}^{+\infty} u'(x)^2\,dx.
$$
More precisely,
\be
\int_{-\infty}^{+\infty} p'(x)^2\,dx \, = \,
\int_{-\infty}^{+\infty} \frac{p'(x)^2}{p(x)}\,p(x)\,dx \, \leq \,
\max_x p(x) \int_{-\infty}^{+\infty} \frac{p'(x)^2}{p(x)}\,dx 
 \, \leq \, I(X)^{3/2}.
\en

Since the estimate on the total variation norm $\|p\|_{{\rm TV}}$
can be given in terms of the Fisher information, it is natural to ask
whether or not it is possible to bound the total variation distance
from $p$ to a normal density in terms of the relative Fisher information.
This suggests the following bound.

\vskip5mm
{\bf Proposition 2.4.} {\it If $X$ has mean zero, variance one, and density $p$ 
with finite Fisher information, then 
\be
\|p - \varphi\|_{{\rm TV}} \leq 4\sqrt{I(X||Z)},
\en
where $Z$ has standard normal density $\varphi$.
}

\vskip5mm
{\bf Proof.} Using
$$
p'(x) - \varphi'(x) = 
\Big(\frac{p'(x)}{p(x)} - \frac{\varphi'(x)}{\varphi(x)}\Big) p(x) -
x\,(p(x) - \varphi(x)) \qquad (p(x)>0)
$$
and applying Cauchy's inequality, we may write
\begin{eqnarray}
\|p - \varphi\|_{{\rm TV}}
 & = &
\int_{-\infty}^{+\infty} |p'(x) - \varphi'(x)|\,dx \nonumber \\
 & \leq &
I(X||Z)^{1/2} +  \int_{-\infty}^{+\infty} |x|\,|p(x) - \varphi(x)|\,dx.
\end{eqnarray}
The last integral represents a weighted total variation distance between the
distributions of $X$ and $Z$ with weight function $w(x) = |x|$. 

On this step we apply the following extention of Csisz\'ar-Kullback-Pinsker's
inequality (CKP) to the scheme of weighted total variation distances, which is
proposed by Bolley and Villani, cf. [B-V], Theorem 2.1 (ii).
If $X$ and $Y$ are random variables with densities $p$ and $q$, and 
$w(x) \geq 0$ is a measurable function, then
$$
\Big(\int_{-\infty}^{+\infty} w(x)\,|p(x) - q(x)|\,dx\Big)^2
 \, \leq \, C D(X||Y) \, = \, 
C \int_{-\infty}^{+\infty} p(x)\,\log\frac{p(x)}{q(x)}\,dx,
$$
where
$$
C \, = \, 2\, \Big(1 + \log \int_{-\infty}^{+\infty} e^{w(x)^2} q(x)\,dx\Big).
$$
The inequality also holds in the setting of abstract measurable spaces,
and when $w=1$ it yields the classical CKP inequality with an additional 
factor $2$. 

In our case, $Y=Z$, $q = \varphi$, and taking 
$w(x) = \sqrt{t/2}\, |x|$ $(0 < t < 1$), we get
$$
\frac{t}{2}\,
\Big(\int_{-\infty}^{+\infty} |x|\,|p(x) - \varphi(x)|\,dx\Big)^2
 \, \leq \, \Big(2 + \log \frac{1}{1-t}\Big)\,D(X||Z).
$$
One may choose, for example, $t = 1 - \frac{1}{e}$, and recalling (1.1), 
we arrive at
$$
\int_{-\infty}^{+\infty} |x|\,|p(x) - \varphi(x)|\,dx \, \leq \,
3.1\,D(X||Z)^{1/2} \leq \frac{3.1}{\sqrt{2}} \, I(X||Z)^{1/2}.
$$
It remains to use this bound in (2.5), and (2.4) follows.


\vskip10mm
\section{{\bf Fisher information as a functional}}
\setcounter{equation}{0}

\vskip2mm
\noindent
It is worthwile to discuss separately a few general properties
of the Fisher information viewed as a functional on the space 
of densities. We start with topological properties.

\vskip5mm
{\bf Proposition 3.1.} {\it Let $(X_n)_{n \geq 1}$ be a sequence of random 
variables, and $X$ be a random variable such that $X_n \Rightarrow X$ weakly 
in distribution. Then
\be
I(X) \leq \liminf_{n \rightarrow \infty} \, I(X_n).
\en
}

\vskip2mm
Denote by ${\mathfrak P}_1$ the collection of all (probability) densities 
on the real line with finite Fisher information, and let ${\mathfrak P}_1(I)$ 
denote the subset of all densities which have Fisher information of at most 
size $I>0$. On the set ${\mathfrak P}_1$ the relation (3.1) may be written as
\be
I(p) \leq \liminf_{n \rightarrow \infty} \, I(p_n),
\en
which holds under the condition that the corresponding distributions are 
convergent weakly, i.e.,
\be
\lim_{n \rightarrow \infty} \, \int_{-\infty}^a p_n(x)\,dx \, = \,
\int_{-\infty}^a p(x)\,dx, \qquad {\rm for \ all} \ \ a \in \R.
\en
Hence, every ${\mathfrak P}_1(I)$ is closed in the weak topology. In fact, 
inside such sets (3.3) can be strengthened to the convergence in the 
$L^1$-metric,
\be
\lim_{n \rightarrow \infty} \,
\int_{-\infty}^{-\infty} |p_n(x)\,dx - p(x)|\,dx \, = \, 0.
\en

\vskip5mm
{\bf Proposition 3.2.} {\it On every set ${\mathfrak P}_1(I)$ the weak topology 
with convergence $(3.3)$ and and the topology generated by the $L^1$-norm 
coincide, and the Fisher information is a lower semi-continuous functional 
on this set.
}

\vskip5mm
{\bf Proof.} For the proof of Proposition 3.1, one may assume that 
$I(X_n) \rightarrow I$, for some (finite) constant $I$. Then, for 
sufficiently large $n$, the $X_n$ have absolutely continuous densities $p_n$ 
with Fisher information at most $I+1$. By Proposition 2.2, such densities 
are uniformly bounded and have uniformly bounded variations. Hence, by the 
second Helly theorem (cf. e.g. [K-F]), there are a subsequence $p_{n_k}$ and 
a function $p$ of bounded variation, such that $p_{n_k}(x) \rightarrow p(x)$, 
as $k \rightarrow \infty$, for all points $x$. Necessarily, $p(x) \geq 0$ 
and $\int_{-\infty}^{+\infty} p(x)\,dx \leq 1$. Since the sequence of 
distributions of $X_n$ is tight (or weakly pre-compact), it also follows that
$\int_{-\infty}^{+\infty} p(x)\,dx = 1$. Hence, $X$ has an absolutely 
continuous distribution with $p$ as its density, and the weak convergence (3.3) 
holds. 

For the proof of Proposition 3.2, a similar argument should be applied to
an arbitrary prescribed subsequence $p_{n_k}$, where we obtain
$p(x) = \lim_{l \rightarrow \infty} p_{n_{k_l}}(x)$ for some further
subsequence. By Scheffe's lemma, this property implies the convergence 
in $L^1$-norm, that is, (3.4) holds along $n_{k_l}$. This implies the 
convergence in $L^1$ for the whole sequence $p_n$, which is the assertion 
of Proposition 3.2.

To continue the proof of Proposition 3.1, for simplicity of notations, 
assume that the subsequence constructed in the first step is actually 
the whole sequence. By (2.3),
$$
\int_{-\infty}^{+\infty} p_n'(x)^2\,dx \leq (I + 1)^{3/2},
$$
which implies that the derivatives are uniformly integrable on every finite 
interval. By the Dunford-Pettis compactness criterion for the space $L^1$
(over finite measures), there is a subsequence $p_{n_k}'$ which is convergent 
to some locally integrable function $u$ in the sense that
\be
\int_A p_{n_k}'(x)\,dx \rightarrow \int_A u(x)\,dx,
\en
for any bounded Borel set $A \subset \R$. (This is the weak 
$\sigma(L^1,L^\infty)$ convergence on finite intervals.)
Note that, according to Proposition 2.1, $p_{n_k}'$ may be replaced in (3.5) 
with the sequence $p_{n_k}' 1_{\{p_{n_k} > 0\}}$, which is thus convergent 
to~$u$ as well. 

Taking a finite interval $A = (a,b)$ in (3.5), we get
$$
\int_a^b u(x)\,dx = p(b) - p(a),
$$
which means that $p$ is (locally) absolutely continuous. Furthermore, since
$$
\|p\|_{{\rm TV}} = \int_{-\infty}^{+\infty} |u(x)|\,dx
$$
is finite, we conclude that $u \in L^1(\R)$, thus representing a Radon-Nikodym
derivative: $u(x) = p'(x)$. Again, for simplicity of notations, assume 
the subsequence of derivatives obtained is actually the whole sequence.

Next, consider the sequence of functions 
$$
\xi_n(x) = \frac{p_n'(x)}{\sqrt{p_n(x)}}\, 1_{\{p_n(x) > 0\}}.
$$
They have $L^2(\R)$-norm bounded by $\sqrt{I+1}$ (for large $n$).
Since the unit ball of $L^2$ is weakly compact, there is a subsequence
$\xi_{n_k}$ which is weakly convergent to some function $\xi \in L^2$, that is,
$$
\int_{-\infty}^{+\infty} \xi_{n_k}(x)\, q(x)\,dx \rightarrow 
\int_{-\infty}^{+\infty} \xi(x)\, q(x)\,dx,
$$
for any $q \in L^2$. As a consequence,
$$
\int_{-\infty}^{+\infty} \xi_{n_k}(x)\, \sqrt{p_{n_k}(x)}\, q(x)\,dx \rightarrow 
\int_{-\infty}^{+\infty} \xi(x)\, \sqrt{p(x)}\, q(x)\,dx,
$$
due to the uniform boundedness and pointwise convergence of $p_n$.
In other words, again omitting sub-indices, the functions 
$p_n'\, 1_{\{p_n > 0\}}$ are weakly convergent in $L^2$ to the function 
$\xi\sqrt{p}$. In particular, for $q = 1_A$ with an arbitrary bounded Borel set 
$A \subset \R$,
$$
\int_A p_n'\, 1_{\{p_n > 0\}}\,dx \rightarrow \int_A \xi(x)\sqrt{p(x)}\ dx.
$$

As a result, we have obtained two limits for $p_n'\, 1_{\{p_n > 0\}}$, 
which must coincide, i.e., we get $\xi\sqrt{p} = u = p'$ a.e. Hence, 
$p = 0 \Rightarrow p' = 0$ and $\xi = \frac{p'}{\sqrt{p}}$ \, a.e. 
on the set $\{p(x) > 0\}$.
Finally, the weak convergence $\xi_{n_k} \rightarrow \xi$ in $L^2$,
as in any Banach space, yields
$$
I(p) \, = \, \|\xi \cdot 1_{\{p>0\}}\|_{L^2}^2 \, \leq \,  
\|\xi\|_{L^2}^2 \, \leq \,
\liminf_{k \rightarrow \infty}\, \|\xi_{n_k}\|_{L^2}^2 \, = \, 
\liminf_{n \rightarrow \infty} \, I(p_{n_k}) \, = \, I.
$$
Thus, Proposition 3.1 is proved.

Another general property of the Fisher information is its convexity, that is,
we have the inequality
\be
I(p) \leq \sum_{i=1}^n \alpha_i I(p_i),
\en
where $p = \sum_{i=1}^n \alpha_i p_i$ with arbitrary densities $p_i$ and 
weights $\alpha_i > 0$, $\sum_{i=1}^n \alpha_i = 1$. This readily follows 
from the fact that the homogeneous function $R(u,v) = u^2/v$ is convex on 
the upper half-plane $u \in \R$, $v>0$. 
Moreover, Cohen [C] showed that the inequality (3.6) is strict.

As a consequence, the collection ${\mathfrak P}_1(I)$ of all densities 
on the real line with Fisher information $\leq I$ represents a convex 
closed set in the space $L^1 = L^1(\R)$ (for strong or weak topologies).

We need to extend Jensen's inequality (3.6) to arbitrary "continuous" convex 
mixtures of densities. In order to formulate this more precisely, recall
the definition of mixtures. Denote by $\mathfrak P$ the collection 
of all densities, which represents a closed subset of $L^1$ with the weak
$\sigma(L^1,L^\infty)$ topology. For any Borel set $A \subset \R$, the 
functionals $q \rightarrow \int_A q(x)\,dx$ are bounded and continuous on 
$\mathfrak P$. So, given a Borel probability measure $\pi$ on ${\mathfrak P}$, 
one may introduce the probability measure on the real line
\be
\mu(A) = \int_{\mathfrak P} \bigg[\int_A q(x)\,dx\bigg]\,d\pi(q).
\en
It is absolutely continuous with respect to Lebesgue measure and has
some density $p(x) = \frac{d\mu(x)}{dx}$ called the (convex) mixture of 
densities with mixing measure $\pi$. For short,
$$
p(x) = \int_{\mathfrak P} q(x)\,d\pi(q).
$$

\vskip2mm
{\bf Proposition 3.3.} {\it If $p$ is a convex mixture of densities 
with mixing measure $\pi$, then
\be
I(p) \leq \int_{\mathfrak P} I(q)\,d\pi(q).
\en
}

\vskip2mm
{\bf Proof.} Note that the integral in (3.8) makes sense, since the
functional $q \rightarrow I(q)$ is lower semi-continuous and hence Borel
measurable on $\mathfrak P$ (Proposition 3.1). We may assume that
this integral is finite, so that $\pi$ is supported on the convex 
(Borel measurable) set $\mathfrak P_1 = \cup_I \mathfrak P_1(I)$. 

Identifying densities with corresponding probability measures (having these 
densities), we consider $\mathfrak P_1$ as a subset of the locally convex 
space $E$ of all finite measures $\mu$ on the real line endowed with the weak 
topology.

{\it Step} 1. Suppose that the measure $\pi$ is supported on some convex 
compact set  $K$ contained in $\mathfrak P_1(I)$. Since the functional 
$q \rightarrow I(q)$ is finite, convex and lower semi-continuous on $K$, 
it admits the representation
$$
I(q) \, = \, \sup_{l \in \mathfrak L} \, l(q), \qquad q \in K,
$$
where $\mathfrak L$ denotes the family of all continuous affine functionals 
$l$ on $E$ such that $l(q) < I(q)$, for all $q \in K$ (cf. e.g. Meyer [M], 
Chapter XI, Theorem T7).
In our particular case, any such functional acts on probability measures as
$l(\mu) = \int_{-\infty}^{+\infty} \psi(x)\,d\mu(x)$
with some bounded continuous function $\psi$ on the real line. Hence,
$$
I(q) \, = \, \sup_{\psi \in \mathfrak C} \, 
\int_{-\infty}^{+\infty} q(x)\psi(x)\,dx,
$$
for some family $\mathfrak C$ of bounded continuous functions $\psi$ on $\R$.
An explicit description of $\mathfrak C$ would be of interest, 
but this question will not be pursued here. As a consequence, 
by the definition (3.7) for the measure $\mu$ with density $p$,
\bee
\int_{\mathfrak P} I(q)\,d\pi(q)
 & \geq &
\sup_{\psi \in \mathfrak C} \ \int_{\mathfrak P} \bigg[
\int_{-\infty}^{+\infty} q(x)\psi(x)\,dx\bigg]\,d\pi(q) \\
 & = &
\sup_{\psi \in \mathfrak C} \ \int_{-\infty}^{+\infty} p(x)\psi(x)\,dx
 \ = \ I(p),
\ene
which is the desired inequality (3.8).

{\it Step} 2. Suppose that $\pi$ is supported on $\mathfrak P_1(I)$, for some 
$I>0$. Since any finite measure on $E$ is Radon, and since the set
$\mathfrak P_1(I)$ is closed and convex, there is an increasing sequence of 
compact subsets $K_n \subset \mathfrak P_1(I)$ such that $\pi(\cup_n K_n) = 1$. Moreover, $K_n$ can be chosen to be convex (since the closure of the convex 
hull will be compact, as well). Let $\pi_n$ denote the normalized restriction 
of $\pi$ to $K_n$ (with sufficiently large $n$ so that $c_n = \pi(K_n) > 0$) 
and define its baricenter
\be
p_n(x) = \int_{K_n} q(x)\,d\pi_n(q).
\en
From (3.7) it follows that the measures with densities $p_n$ are weakly
convergent to the measure $\mu$ with density $p$, hence the relation (3.2) 
holds: $I(p) \leq \liminf_{n \rightarrow \infty} \, I(p_n)$.
On the other hand, by the previous step,
\be
I(p_n) \, \leq \, \int_{K_n} I(q)\,d\pi_n(q) = \frac{1}{c_n}\, 
\int_{K_n} I(q)\,d\pi(q) \, \rightarrow \int_{\mathfrak P_1(I)} I(q)\,d\pi(q),
\en
which yields (3.8).

{\it Step} 3. In the general case, we may apply Step 2 to 
the normalized restrictions $\pi_n$ of $\pi$ to the sets $K_n = \mathfrak P_1(n)$.
Again, for the densities $p_n$ defined as in (3.9), we obtain (3.10),
where $\mathfrak P_1(I)$ should be replaced with $\mathfrak P_1$.
Another application of the lower semi-continuity of the Fisher information
finishes the proof.


\vskip10mm
\section{{\bf Convolution of three densities of bounded variation}}
\setcounter{equation}{0}

\vskip2mm
\noindent
Although densities with finite Fisher information must be functions
of bounded variation, the converse is not always true. Nevertheless,
starting from a density of bounded variation and taking several 
convolutions with itself, the resulting density will have finite 
Fisher information. Our nearest aim is to prove:

\vskip5mm
{\bf Proposition 4.1.} {\it If independent random variables $X_1,X_2,X_3$
have densities $p_1,p_2,p_3$ with finite total variation, then $S = X_1 + X_2 + X_3$
has finite Fisher information, and moreover,
\be
I(S) \, \leq \, \frac{1}{2}\,\Big[
\|p_1\|_{{\rm TV}} \, \|p_2\|_{{\rm TV}} +
\|p_1\|_{{\rm TV}} \, \|p_3\|_{{\rm TV}} +
\|p_2\|_{{\rm TV}} \, \|p_3\|_{{\rm TV}}\Big].
\en
}

\vskip2mm
One may further extend (4.1) to sums of more than 3 independent summands,
but this will not be needed for our purposes (since the Fisher information
may only decrease when adding an independent summand.)

In the i.i.d. case the above estimate can be simplified. By a direct 
application of the inverse Fourier formula, the right-hand side of (4.1) may 
be related furthermore to the characteristic functions of $X_j$. 
We will return to this in the next section.

First let us look at the particular case where $X_j$ are uniformly distributed
over intervals. This important example already shows that the Fisher
information $I(X_1 + X_2)$ does not need to be finite, while it is finite 
for 3 summands. (This somewhat curious fact was pointed out to one of the 
authors by K. Ball.) In fact, there is a simple quantitative bound.

\vskip5mm
{\bf Lemma 4.2.} {\it If independent random variables $X_1,X_2,X_3$
are uniformly distributed on intervals of lengths $a_1,a_2,a_3$, then
\be
I(X_1 + X_2 + X_3) \, \leq \, 2\, 
\bigg[\frac{1}{a_1 a_2} + \frac{1}{a_1 a_3} + \frac{1}{a_2 a_3}\bigg].
\en
}

\vskip2mm
The density of the sum $S = X_1 + X_2 + X_3$ may easily be evaluated and 
leads to a rather routine problem of estimation of 
$I(S)$ as a function of the parameters $a_j$. Alternatively, there is 
an elegant approach based on general properties of so-called convex or 
hyperbolic distributions and the fact that the density $p$ of $S$ behaves 
like the beta density near the end points of the supporting interval. 

To describe the argument, let us recall a few definitions and results 
concerning such measures. A probability measure $\mu$ on $\R^d$ is called
$\kappa$-concave with a (convexity) parameter $0 < \kappa \leq 1$, 
if it satisfies a Brunn-Minkowski-type inequality
$$
\mu(tA + (1-t)B) \geq \big(t \mu(A)^\kappa + (1-t) \mu(B)^\kappa\big)^{1/\kappa}
$$
in the class of all non-empty Borel sets $A,B \subset \R^d$, and for
arbitrary $0 < t < 1$. We refer to the papers by Borell [Bor1-2] for 
basic properties of such measures, cf. also [Bo] (in fact, the values 
$\kappa \leq 0$ are also allowed, but will not be needed here).

If $\mu$ is absolutely continuous, the definition reduces to the property
that $\mu$ is supported on some open convex set $\Omega \subset \R^d$
(necessarily bounded), where it has a positive density $p$ such that the 
function $p^{\kappa/(1-\kappa d)}$ is concave on $\Omega$ (Borell's 
characterization theorem). For example, the normalized Lebesgue measure on 
any convex body is $\frac{1}{d}$-concave. In dimension one, $\mu$ has to be 
supported on some finite interval $(x_0,x_1)$, and Borell's description may 
also be given in terms of the function
$$
L(t) = p(F^{-1}(t)), \qquad 0 < t < 1,
$$
where $F^{-1}:(0,1) \rightarrow (x_0,x_1)$ denotes the inverse of the 
distribution function $F(x) = \mu(x_0,x)$, restricted to the supporting 
interval. Namely (cf. [Bo]), a probability measure $\mu$ is $\kappa$-concave, 
if and only if the function $L^{1/(1-\kappa)}$ is concave on $(0,1)$.

We only need the following well-known fact about the convexity parameter
of convolutions which we formulate in case of three measures: 
If $\mu_j$ are $\kappa_j$-concave $(j = 1,2,3)$, then
the measure $\mu = \mu_1 * \mu_2 * \mu_3$ is $\kappa$-concave, where
\be
\frac{1}{\kappa} = \frac{1}{\kappa_1} + \frac{1}{\kappa_2} + \frac{1}{\kappa_3}.
\en

Note also that the Fisher information of a random variable $X$ with density $p$ 
is expressed in terms of the associated function $L$ as
\be
I(X) = \int_0^1 L'(t)^2\,dt.
\en
This general formula holds whenever $p$ is absolutely continuous and positive
on the supporting interval (without any $\kappa$-concavity assumption).

\vskip5mm
{\bf Proof of Lemma 4.2.} For definiteness, let $X_j$ take values in $[0,a_j]$. 
Since the distributions of $X_j$ are $1$-concave, the distribution of
$S = X_1 + X_2 + X_3$ is $\frac{1}{3}$-concave, according to (4.3).
This means that $S$ has density $p$ such that $p^{1/2}$ is concave on the
supporting interval $[0,a_1+a_2+a_3]$, or equivalently, $L^{3/2}$ is concave 
on $(0,1)$, where $L$ is the associated function for $S$.

Note that $S$ has an absolutely continuous density $p$, which is thus vanishing 
at the end points $x = 0$ and $x = a_1 + a_2 + a_3$. Hence, $L(0+) = L(1-) = 0$. 
By the concavity, the Radon-Nikodym derivative
$(L^{3/2})' = \frac{3}{2}\,L^{1/2}\, L'$ is non-increasing, and since
$L$ is symmetric about the point $\frac{1}{2}$, we get, for all $0 < t < 1$,
$$
L'(t)^2\, L(t) \, \leq \, c, \qquad {\rm where} \ \ \ \
c \, = \, \lim_{t \rightarrow 0}\, L'(t)^2\, L(t).
$$
Hence, by (4.4), 
\be
I(X) \, \leq \, \int_0^1 \frac{c}{L(t)}\,dt \, = \, c\,(a_1 + a_2 + a_3).
\en

It remains to find the constant $c$. Putting $a = a_1 a_2 a_3$, it should be
clear that, for all $x > 0$ and $t>0$ small enough,
$$
F(x) = \P\{S \leq x\} = \frac{x^3}{6a}, \quad p(x) = \frac{x^2}{2a}, 
\quad F^{-1}(t) = (6at)^{1/3}, \quad L(t) = \frac{1}{2a}\,(6at)^{2/3},
$$
and finally $c = L'(t)^2\, L(t) = \frac{2}{a}$. Thus, in (4.5) we arrive at
$I(X) \leq \frac{2}{a}\,(a_1 + a_2 + a_3)$ which is exactly (4.2).

\vskip2mm
Lemma 4.2 allows us to reduce Proposition 4.1 to the case of uniform
distrubutions. Note that if a density $p$ is written as a convex mixture
\be
p(x) = \int_{\mathfrak P} q(x)\,d\pi(q),
\en
then by the convexity of the total variation norm,
\be
\|p\|_{\rm TV} \leq \int_{\mathfrak P} \|q\|_{\rm TV}\,d\pi(q).
\en
Recall that we understand (4.6) as the equality (3.7) of the corresponding 
measures. So, (4.7) is also uses our original agreement that, for each $x$, 
the value $p(x)$ lies in the closed segment with endpoints $p(x-)$ and $p(x+)$.

In order to apply Lemma 4.2 together with Jensen's inequality for Fisher 
information, we need however to require that $\pi$ has to be supported on 
uniform densities (that is, densities of normalized Lebesgue measures on 
finite intervals) and secondly to reverse (4.7). Indeed this turns out 
to be possible, which may be a rather interesting observation.

\vskip5mm
{\bf Lemma 4.3.} {\it Any density $p$ of bounded variation 
can be represented as a convex mixture $(4.6)$ of uniform densities with 
a mixing measure $\pi$ such that
\be
\|p\|_{\rm TV} = \int_{\mathfrak P} \|q\|_{\rm TV}\,d\pi(q).
\en
}

\vskip2mm
For example, if $p$ is supported and non-increasing on $(0,+\infty)$,
there is a canonical representation
$$
p(x) = \int_0^{+\infty} \frac{1}{x_1}\,1_{\{0 < x < x_1\}}\,d\pi(x_1) 
\qquad {\rm a.e.}
$$
with a unique mixing probability measure $\pi$ on $(0,+\infty)$. In this case
$\|p\|_{\rm TV} = 2p(0+)$, and (4.8) is obvious. One may write
a similar representation for densities of unimodal distributions.
In general, another way to write (4.6) and (4.8) is
\bee
p(x) & = & 
\int_{x_1>x_0} \frac{1}{x_1 - x_0}\,1_{\{x_0 < x < x_1\}}\, d\pi(x_0,x_1), \\
\|p\|_{\rm TV}
     & = &
2 \int_{x_1>x_0} \frac{1}{x_1 - x_0}\ d\pi(x_0,x_1),
\ene
where $\pi$ is a 
Borel probability measure on the half-plane $x_1 > x_0$ 
(i.e., above the main diagonal).

Let us also note that the sets ${\rm BV}(c)$ of all densities $p$ with 
$\|p\|_{\rm TV} \leq c$ are closed under the weak convergence (3.3) of the 
corresponding probability distributions. Moreover, the weak convergence in
${\rm BV}(c)$ coincides with convergence in $L^1$-norm, which can be proved 
using the same arguments as in the proof of Proposition 3.2. In particular, the 
functional $q \rightarrow \|q\|_{\rm TV}$ is lower semi-continuous and hence 
Borel measurable on $\mathfrak P$, so the integrals (4.7)-(4.8) make sense. 

Denote by $U$ the collection of all uniform densities which thus may be 
identified with the half-plane $\tilde U = \{(a,b) \in \R^2: b > a\}$ via 
the map $(a,b) \rightarrow q_{a,b}(x) = \frac{1}{b-a}\,1_{\{a<x<b\}}$.
The usual convergence on $\tilde U$ in the Euclidean metric coincides 
with the weak convergence (3.3) of $q_{a,b}$. The closure of $U$ for the 
weak topology contains $U$ and all delta-measures, hence $U$ is a Borel 
measurable subset of $\mathfrak P$.

\vskip2mm
{\bf Proof.} We only need the existence part which is proved below in two steps.

{\it Step} 1. First consider the discrete case, where $p$ is piecewise 
constant, i.e., it is supported and constant on consecutive semiopen intervals 
$\Delta_k = [x_{k-1},x_k)$, $k = 1,\dots, n$, where $x_0 < ... < x_n$. 
Putting $p(x) = c_k$ on $\Delta_k$, we then have
$$
\|p\|_{\rm TV} = c_1 + |c_2 - c_1| + \dots + |c_n - c_{n-1}| + c_n.
$$

In this case the existence of the representation (4.6), moreover -- with 
a discrete mixing measure $\pi$, satisfying (4.8), can be proved by induction 
on $n$. If $n = 1$ or $n = 2$, then $p$ is monotone on $\Delta_1$, 
respectively, on $\Delta_1 \cup \Delta_2$, and the statement is obvious. 

If $n \geq 3$, one should distinguish between several cases.
If $c_1 = 0$ or $c_n = 0$, we are reduced to the smaller number of
supporting intervals. If $c_k = 0$ for some $1 < k < n$, one can write
$p = f + g$ with $f(x) = p(x)\,1_{\{x<x_{k-1}\}}$, 
$g(x) = p(x)\,1_{\{x \geq x_k\}}$. These functions are supported on
disjoint half-axes, so $\|p\|_{\rm TV} = \|f\|_{\rm TV} + \|g\|_{\rm TV}$.
Moreover, the induction hypothesis may be applied to both $f$ and $g$
(or one can first normalize these functions to work with densities, but
this is less convenient). As a result,
$$
f = f_1 + \dots + f_k, \qquad g = g_1 + \dots + g_l \quad {\rm a.e.}
$$
where each $f_i$ is supported and constant on some interval inside $[x_0,x_{k-1})$,
each $g_j$ is supported and constant on some interval inside $[x_k,x_n)$,
and 
$$
\|f\|_{\rm TV} = \|f_1\|_{\rm TV} + \dots + \|f_k\|_{\rm TV}, \qquad 
\|g\|_{\rm TV} = \|g_1\|_{\rm TV} + \dots + \|g_l\|_{\rm TV} .
$$
Hence, 
$$
p = \sum_i f_i + \sum_j g_j \quad {\rm with} \quad
\|f\|_{\rm TV} = \sum_i \|f_i\|_{\rm TV} + \sum_j \|g_j\|_{\rm TV}.
$$

Finally, assume that $c_k > 0$ for all $k \leq n$. Putting $c_* = \min_k c_k$,
write $p = f + g$, where $f = c_*\, 1_{[x_0,x_n)}$ and 
$g$ thus takes the values $c_k - c_*$ on $\Delta_k$. Clearly,
$$
\|p\|_{\rm TV} = 2c_* + \|g\|_{\rm TV} = \|f\|_{\rm TV} + \|g\|_{\rm TV}.
$$
By the definition, $g$ takes the value zero on one of the intervals
(where $c_k = c_*$), so we are reduced to the previous step. On that step, 
we obtained a representation $g = g_1 + \dots + g_l$ such that 
$\|g\|_{\rm TV} = \|g_1\|_{\rm TV} + \dots + \|g_l\|_{\rm TV}$, where each 
$g_j$ is supported and constant on some interval inside $[x_0,x_n)$. Hence,
$$
p = f + \sum_j g_j \quad {\rm with} \quad
\|p\|_{\rm TV} = \|f\|_{\rm TV} + \sum_j \|g_j\|_{\rm TV}.
$$

Although the measure $\pi$ has not been constructed constructively,
one may notice that it should be supported on the densities of the form
$$
q_{ij}(x) = \frac{1}{x_j - x_i}\,1_{\{x_i \leq x < x_j\}}, \qquad
0 \leq i < j \leq n.
$$

{\it Step} 2. In the general case, one may assume that $p$ is right-continuous.
Consider the collection of piecewise constant densities of the form
\be
\tilde p(x) = d \, \sum_{k=1}^n p(x_{k-1}) \,1_{\{x_{k-1} \leq x < x_k\}}
\en
with arbitrary points $x_0 < ... < x_n$ of continuity of $p$ 
such that $p(x_{k-1}) > 0$ for at least one $k$, and where $d$ is 
a normalizing constant so that $\int_{-\infty}^{+\infty} \tilde p(x)\,dx = 1$.
Since $p$ has bounded total variation, it is possible to
construct a sequence $p_n$ of the form $(4.9)$ which is convergent to $p$
in $L^1$-norm and with $d = d_n \rightarrow 1$. By the construction,
\be
\frac{1}{d_n}\,\|p_n\|_{\rm TV} \, = \, p(x_0) + p(x_{n-1}) + 
\sum_{k=1}^{n-1} |p(x_k) - p(x_{k-1})| \, \leq \, \|p\|_{\rm TV},
\en
so all $p_n$ belong to ${\rm BV}(c)$ with some constant $c$.

Using the previous step, one can define discrete probability 
measures $\pi_n$ supported on $U$ and such that
\be
p_n(x) = \int_U q(x)\,d\pi_n(q), \qquad
\|p_n\|_{\rm TV} = \int_U \|q\|_{\rm TV}\,d\pi_n(q).
\en
Since $U$ has been identified with the half-plane $\tilde U$,
replacing $d\pi_n(q)$ with $d\pi_n(a,b)$ should not lead to confusion.
In particular, the second equality in (4.11) may be written as
\be
\|p_n\|_{\rm TV} = 2\int_{\tilde U} \frac{1}{b-a}\,d\pi_n(a,b).
\en

From the first equality in (4.11) it follows that, for any $T>0$,
$$
\int_U \Big[\int_{|x| \geq T} q(x)\,dx\Big]\,d\pi_n(q) \, = \,
\int_{|x| \geq T} p_n(x) \, \leq \, \int_{|x| \geq T} p(x)\,dx + \|p_n - p\|_1.
$$
Hence, by Chebyshev's inequality, for any $\ep_k > 0$,
\be
\pi_n\Big\{q \in U: \int_{|x| \geq k} q(x)\,dx > \ep_k\Big\} \, \leq \,
\frac{1}{\ep_k}\,\Big(\int_{|x| \geq k} p(x)\,dx + \|p_n - p\|_1\Big).
\en
Clearly, one can choose a sequence $\ep_k \downarrow 0$ and an increasing 
sequence of indices $n_k$ such that the right-hand side of (4.13) will tend 
to zero, as $k \rightarrow \infty$, uniformly over all $n \geq n_k$. 
In particular, the above inequality holds for $\pi_{n_k}$.

On the other hand (identifying $q$ with corresponding probability distributions), 
by the Prokhorov compactness criterion, the collection of densities 
$$
\Big\{q \in {\mathfrak P}: \int_{|x| \geq k} q(x)\,dx \leq \ep_k\Big\}
$$
is pre-compact for the weak topology with convergence (3.3), cf. e.g. [Bi].
Therefore, by the same criterion applied to $\mathfrak P$ as a Polish space, 
$\pi_n$ contains a weakly convergent subsequence $\pi_{n_k}$
with some limit $\pi \in \mathfrak P$. This measure is supported on the (weak) 
closure of $U$, which is a larger set, since it contains delta-measures, or 
the main diagonal in $\R^2$, if we identify $U$ with $\tilde U$. However, 
using (4.12) together with Chebyshev's inequality, and then applying (4.10), 
we see that, for any $\ep > 0$ and all $n \geq n_0$,
$$
\pi_n\{(a,b): b - a < \ep\} \, = \,
\pi_n\Big\{(a,b): \frac{1}{b - a} > \frac{1}{\ep}\Big\} \, \leq \,
\frac{\ep}{2}\, \|p_n\|_{\rm TV} \, < \, \ep\,\|p\|_{\rm TV}.
$$
Hence, $\pi$ is actually supported on $U$. Moreover, taking the limit along 
$n_k$ in the first equality in (4.11), we obtain the representation (4.6).

Now, the sets $G(t) = \{q \in U:\|q\|_{\rm TV} > t\}$ are open in the weak 
topology (by the lower semicontinuity of the total variation norm), hence,
$\liminf_{k \to \infty} \pi_{n_k}(G(t)) \geq \pi(G(t))$. Applying
Fatou's lemma and then again (4.10) and the second equality in (4.11), we get
\bee
\int_U \|q\|_{\rm TV}\,d\pi(q)
 & = &
\int_0^{+\infty} \pi(G(t))\,dt  \ \leq \ \liminf_{k \to \infty}\, 
\int_0^{+\infty} \pi_{n_k}(G(t))\,dt \\
 & = &
\liminf_{k \to \infty}\, \int_U \|q\|_{\rm TV}\,d\pi_{n_k}(q)
 \ = \ 
\liminf_{k \to \infty}\, \|p_{n_k}\|_{\rm TV} \ \leq \ \|p\|_{\rm TV}.
\ene
In view of Jensen's inequality (4.7), we obtain (4.8) thus proving the existence 
part of the lemma.

\vskip2mm
{\bf Proof of Proposition 4.1.} 
We may write down the representation (4.6) from Lemma 4.2 for each of the 
densities $p_j$ $(j=1,2,3)$. That is,
$$
p_j(x) = \int q(x)\,d\pi_j(q) \qquad {\rm a.e.}
$$
with some mixing probability measures $\pi_j$, supported on $U$ and satisfying
\be
\|p_j\|_{\rm TV} = \int \|q\|_{\rm TV}\,d\pi_j(q).
\en 
Taking the convolution, we then have a similar representation
$$
(p_1 * p_2 * p_3)(x) \, = \, \int \!\!\int \!\!\int
(q_1 * q_1 * q_3)(x)\ d\pi_1(q_1) d\pi_2(q_2) d\pi_3(q_3)
\quad {\rm a.e.}
$$
One can now use Jensen's inequality (3.8) for the Fisher information 
and apply (4.2) to bound $I(p_1 * p_2 * p_3)$ from above by
$$
\frac{1}{2} \int \!\!\int \!\!\int \big[
\|q_1\|_{{\rm TV}} \, \|q_2\|_{{\rm TV}} +
\|q_1\|_{{\rm TV}} \, \|q_3\|_{{\rm TV}} +
\|q_2\|_{{\rm TV}} \, \|q_3\|_{{\rm TV}}\big]
\ d\pi_1(q_1) d\pi_2(q_2) d\pi_3(q_3).
$$
In view of (4.14), the triple integral coincides with the right-hand of (4.1).

Proposition 4.1 is proved.


\vskip10mm
\section{{\bf Bounds in terms of characteristic functions}}
\setcounter{equation}{0}

\vskip2mm
\noindent
In view of Proposition 4.1, let us describe how to bound the total variation 
norm of a given density $p$ of a random variable $X$ in terms of the 
characteristic function $f(t) = \E\, e^{itX}$. There are many different bounds
depending on the integrability properties of $f$ and its derivatives, which 
may also depend on assumptions on the finiteness of moments of $X$. We shall 
present two of them here. 

Recall that, if $p$ is absolutely continuous, then
$$
\|p\|_{{\rm TV}} \, = \, \int_{-\infty}^{+\infty} |p'(x)|\,dx.
$$

\vskip5mm
{\bf Proposition 5.1.} {\it If $X$ has finite second moment and
\be
\int_{-\infty}^{+\infty} |t|\,\big(|f(t)| + |f'(t)| + |f''(t)|\big)\,dt < +\infty,
\en
then $X$ has a continuously differentiable density $p$ with finite 
total variation
\be
\|p\|_{{\rm TV}} \, \leq \, \frac{1}{2}\, 
\int_{-\infty}^{+\infty} \big(|tf''(t)| + 2\,|f'(t)| + |t f(t)|\big)\,dt.
\en
}

\vskip2mm
{\bf Proof.} The argument is standard, and we recall it here for completeness.
 
First, by the moment assumption, $f$ is twice continuously differentiable.
The assumption (5.1) implies that $X$ has a continuously differentiable 
density
\be
p(x) = \frac{1}{2\pi}\, \int_{-\infty}^{+\infty} e^{-itx} f(t)\,dt
\en
with derivative
\be
p'(x) = -\frac{i}{2\pi}\, \int_{-\infty}^{+\infty} e^{-itx}\, tf(t)\,dt.
\en

Necessarily $f(t) \rightarrow 0$, as $|t| \rightarrow +\infty$, and
the same is true for $f'(t)$ and $f''(t)$. Therefore, one may integrate 
in (5.3) by parts to get, for all $x \in \R$,
\be
x p(x) = - \frac{i}{2\pi}\, \int_{-\infty}^{+\infty} e^{-itx} f'(t)\,dt
\en
and
$$
x^2 p(x) = -\frac{1}{2\pi}\, \int_{-\infty}^{+\infty} e^{-itx} f''(t)\,dt.
$$
By (5.1), we are allowed to differentiate the last equality by performing
differentiation under the integral sign, which together with (5.4) and (5.5) 
gives
$$
(1+x^2) p'(x) \, = \, \frac{i}{2\pi}\, 
\int_{-\infty}^{+\infty} e^{-itx}\, \big(tf''(t) + 2f'(t) - tf(t)\big)\,dt.
$$
Hence, $|p'(x)| \leq \frac{C}{2\pi\,(1 + x^2)}$ with a constant
described as the integral in (5.2). After integration of this pointwise 
bound, the proposition follows.

\vskip2mm
One can get rid of the assumption of existing second derivative in the bound 
above and remove any moment assumption in Proposition 5.1. But we still need
to insist on the corresponding integrability requirements for the 
characteristic function including its differentiability on the positive half-axis.

\vskip5mm
{\bf Proposition 5.2.} {\it Assume the characteristic function $f(t)$ of a random 
variable $X$ has a continuous derivative for $t>0$, with
\be
\int_{-\infty}^{+\infty} t^2\,\big(|f(t)|^2 + |f'(t)|^2\big)\,dt < +\infty.
\en
Then $X$ has an absolutely continuous distribution with density $p$ of bounded
total variation such that
\be
\|p\|_{{\rm TV}} \, \leq \, \bigg(\int_{-\infty}^{+\infty} |t f(t)|^2\,dt
\int_{-\infty}^{+\infty} |(tf(t))'|^2\,dt\bigg)^{1/4}.
\en
}

\vskip2mm
{\bf Proof.} First assume additionally that $f$ and $f'$ decay at infinity 
sufficiently fast (so that $tf(t) \rightarrow 0$, as $|t| \rightarrow +\infty$).
Integrating by parts in (5.4) and since $(t f(t))'$
is integrable near zero, we get a similar representation
$$
x p'(x) = -\frac{1}{2\pi} \int_{-\infty}^{+\infty} e^{-itx}\, (t f(t))'\,dt.
$$
As usual, write $|p'(x)| = \frac{1}{|1 + ix|}\,|(1 + ix) p(x)|$
and use Cauchy's inequality together with Plancherel's formula, to get
\bee
\bigg(\int_{-\infty}^{+\infty} |p'(x)|\,dx\bigg)^2
 & \leq &
\int_{-\infty}^{+\infty} \frac{dx}{1+x^2} \
\int_{-\infty}^{+\infty} (1+x^2)\,p'(x)^2\,dx \\
 & = &
\frac{1}{2}\,\int_{-\infty}^{+\infty} \big[|t f(t)|^2 + |(tf(t))'|^2 \big]\,dt.
\ene
Applying the same inequality to $\lambda X$ and optimizing over $\lambda > 0$,
we arrive at (5.7).

In the general case, one may apply (5.7) to the regularized random
variables $X_\sigma = X + \sigma Z$ with small parameters $\sigma>0$, where
$Z \sim N(0,1)$ is independent of $X$. They have smooth densities $p_\sigma$
and characteristic functions $f_\sigma(t) = f(t)\, e^{-\sigma^2 t^2/2}$.
Repeating the previous argument for the difference of densities, we obtain
an analogue of (5.7),
\be
\|p_{\sigma_1} - p_{\sigma_2}\|_{{\rm TV}}^4 \, \leq \, 
\int_{-\infty}^{+\infty} |t\, (f_{\sigma_1}(t) - f_{\sigma_2}(t))|^2\,dt
\int_{-\infty}^{+\infty} |(t\, (f_{\sigma_1}(t) - f_{\sigma_2}(t)) )'|^2\,dt
\en
with arbitrary $\sigma_1,\sigma_2 > 0$. Since the integrals in (5.7) are 
finite, by the Lebesgue dominated convergence theorem, the right-hand side 
of (5.8) tends to zero, as long as $\sigma_1,\sigma_2 \rightarrow 0$. Hence, 
the family $\{p_\sigma\}$ is fundamental (Cauchy) for 
$\sigma \rightarrow 0$ in the Banach space of all functions of bounded 
variation on the real line that are vanishing at infinity. As a result, 
there exists the limit $p = \lim_{\sigma \rightarrow 0} p_\sigma$ 
in this space in total variation norm. 

Necessarily, $p(x) \geq 0$ for all $x$, and  
$\int_{-\infty}^{+\infty} p(x)\,dx = 1$.
Hence, $X$ has an absolutely continuous distribution with density $p$.
In addition, by (5.7) applied to $p_\sigma$,
$$
\|p\|_{{\rm TV}} \, = \, \lim_{\sigma \rightarrow 0} \, \|p_\sigma\|_{{\rm TV}}
\, \leq \, \lim_{\sigma \rightarrow 0} \,
\bigg(\int_{-\infty}^{+\infty} |t f_\sigma(t)|^2\,dt
\int_{-\infty}^{+\infty} |(tf_\sigma(t))'|^2\,dt\bigg)^{1/4}.
$$
The last limit exists and coincides with the right-hand side of (5.7).

\vskip5mm
{\bf Corollary 5.3.} {\it If the independent random variables $X_1,X_2,X_3$
have finite first absolute moment and a common characteristic function $f(t)$, 
then 
$$
I(X_1 + X_2 + X_3) \, \leq \, \frac{3}{2}\, 
\bigg(\int_{-\infty}^{+\infty} |t f(t)|^2\,dt
\int_{-\infty}^{+\infty} |(tf(t))'|^2\,dt\bigg)^{1/2}.
$$
If $X_1$ has finite second moment, we also have
$$
I(X_1 + X_2 + X_3) \, \leq \, \frac{3}{8}\,\bigg(\int_{-\infty}^{+\infty} 
\big(|tf''(t)| + 2\,|f'(t)| + |t f(t)|\big)\,dt\bigg)^2.
$$
}


\section{{\bf Classes of densities representable as convolutions}}
\setcounter{equation}{0}

\vskip2mm
\noindent
General bounds like those in Proposition 2.1
may considerably be sharpened in the case where $p$ is representable 
as convolution of several densities with finite Fisher information.

\vskip5mm
{\bf Definition 6.1.} Given an integer $k \geq 1$ and a real number $I > 0$, 
denote by ${\mathfrak P}_k(I)$ the collection of all functions $p$ 
on the real line which can be represented as convolution of $k$ 
probability densities with Fisher information at most $I$.

\vskip5mm
Correspondingly, let ${\mathfrak P}_k$ denote the collection of all functions 
$p$ representable as convolution of $k$ probability densities with finite 
Fisher information. 

The collection ${\mathfrak P}_1$ of all densities with finite Fisher 
information has been already discussed in connection with general properties 
of the functional $I$. For growing $k$, the classes ${\mathfrak P}_k(I)$ 
decrease, since the Fisher information may only decrease when adding an 
independent summand. This also follows from the following general inequality 
of Stam
\be
\frac{1}{I(X+Y)} \geq \frac{1}{I(X)} + \frac{1}{I(Y)},
\en
which holds for all independent random variables (cf. [St], [Bl], [J]).
Moreover, it implies that $p = p_1 * \dots * p_k \in {\mathfrak P}_k(I/k)$,
as long as $p_i \in {\mathfrak P}_1(I)$, $i = 1,\dots,k$.

Any function $p$ in ${\mathfrak P}_k$ is $k-1$ times differentiable, and 
its $(k-1)$-th derivative is absolutely continuous and has a Radon-Nikodym 
derivative, denoted by $p^{(k)}$.
Let us illustrate this property in the important case $k=2$. Write
\be
p(x) = \int_{-\infty}^{+\infty} p_1(x-y) p_2(y)\,dx
\en
in terms of absolutely continuous densities $p_1$ and $p_2$ of independent 
summands $X_1$ and $X_2$ of a random variable $X$ with density $p$. 
Differentiating under the integral sign,
we obtain a Radon-Nikodym derivative of the function $p$,
\be
p'(x) = \int_{-\infty}^{+\infty} p_1'(x-y) p_2(y)\,dy =
\int_{-\infty}^{+\infty} p_1'(y) p_2(x-y)\,dy.
\en
The latter expression shows that $p'$ is absolutely continuous 
and has a Radon-Nikodym derivative
\be
p''(x) = \int_{-\infty}^{+\infty} p_1'(y) p_2'(x-y)\,dy,
\en
which is well-defined for all $x$.
In other words, $p''$ appears as the convolution of the functions $p_1'$ 
and $p_2'$ (which are integrable, according to Proposition 2.2).

These formulas may be used to derive a number of elementary relations within 
the class ${\mathfrak P}_k$, and here we shall describe some of them 
for the cases ${\mathfrak P}_2$ and ${\mathfrak P}_3$.

\vskip5mm
{\bf Proposition 6.2.} {\it Given a density $p \in {\mathfrak P}_2(I)$, for all 
$x \in \R$,
\be
|p'(x)| \leq \,I^{3/4} \sqrt{p(x)} \leq I.
\en
Moreover, $p'$ has finite total variation
$$
\|p'\|_{\rm TV} = \int_{-\infty}^{+\infty} |p''(x)|\,dx \leq I.
$$
}

\vskip5mm
The last bound immediately follows from (6.4) and Proposition 2.2.
To obtain the pointwise bound on the derivative, we may appeal to 
Proposition 2.1 and rewrite the first equality in (6.3) as
$$
p'(x) = \int_{-\infty}^{+\infty} 
\frac{p_1'(x-y)}{\sqrt{p_1(x-y)}}\ 1_{\{p_1(x-y) > 0\}} \
\sqrt{p_1(x-y)}\, p_2(y)\,dy.
$$
Using Cauchy's inequality, we get
\bee
p'(x)^2
 & \leq &
I(X_1) \int_{-\infty}^{+\infty} p_1(x-y)\, p_2(y)^2\,dy \\
 & \leq &
I(X_1)\, \max_y p_2(y)\, \int_{-\infty}^{+\infty} p_1(x-y)\, p_2(y)\,dy
 \ \leq \ I(X_1) I(X_2)^{1/2}\, p(x),
\ene
where we applied Proposition 2.2 to the random variable $X_2$ on the last 
step. This gives the first inequality in (6.5), while the second follows
from $p(x) \leq \sqrt{I}$.

Now, we state similar bounds for the second derivative.

\vskip5mm
{\bf Proposition 6.3.} {\it For any density $p \in {\mathfrak P}_2(I)$, 
we have $p(x) = 0 \Rightarrow p''(x) = 0$ and 
$|p''(x)| \leq I^{3/2}$, for all $x$.
In addition,
$$
\int_{\{p(x)>0\}} \frac{p''(x)^2}{p(x)}\,dx \leq I^2.
$$
}

\vskip2mm
{\bf Proof.} Let us start with the representation (6.4) for a fixed value
$x \in \R$. Note that the function $p_1'(x-y)\, p_2'(y)$ appearing in this 
formula is continuous in $y$. By Proposition 2.1, the integral in (6.4) 
may be restricted to the set $\{y:p_2(y)>0\}$. By the same reason, it may also 
be restricted to the set $\{y:p_1(x-y)>0\}$. Hence,
\be
p''(x) = \int_{-\infty}^{+\infty} p_1'(y) p_2'(x-y)\,1_A(y)\,dy,
\en
where $\{y: p_1(x-y)p_2(y)>0\}$. On the other hand, by the
definition (6.2), the assumption $p(x) = 0$ implies that
$p_1(y) p_2(x-y) = 0$ for almost all $y$. Therefore, $1_A(y) = 0$ a.e.,
and thus the integral (6.6) is vanishing, that is, $p''(x) = 0$.

Using the representation (6.4), the bound $|p''(x)| \leq I^{3/2}$ follows from 
the uniform bound (6.5) on $p'$ and the integral bound of Proposition 2.2.

Next, introduce the functions
$u_i(x) = \frac{p_i'(x)}{\sqrt{p_i(x)}}\, 1_{\{p_i(x) > 0\}}$ ($i = 1,2$)
and rewrite (6.4) as
$$
p''(x) = \int_{-\infty}^{+\infty} \big(u_1(x-y) u_2(y)\big) \,
\sqrt{p_1(x-y)p_2(y)} \ dy.
$$
By Cauchy's inequality,
\be
p''(x)^2 \, \leq \,
\int_{-\infty}^{+\infty} u_1(x-y)^2\, u_2(y)^2\, dy \,
\int_{-\infty}^{+\infty} p_1(x-y)p_2(y) \, dx \, = \, u(x)^2 p(x),
\en
where we used $u \geq 0$ given by
\be
u(x)^2 = \int_{-\infty}^{+\infty} u_1(x-y)^2\, u_2(y)^2\, dy.
\en
Clearly,
$$
\int_{-\infty}^{+\infty} u(x)^2\, dx = I(X_1) I(X_2) \leq I^2,
$$
which is the inequality of the proposition.

\vskip5mm
{\bf Proposition 6.4.} {\it Given a density $p \in {\mathfrak P}_3(I)$, we have,
for all $x$,
$$
|p''(x)| \leq I^{5/4} \sqrt{p(x)}.
$$
}

\vskip2mm
Indeed, by the assumption, one may write $p = p_1 * p_2$ with 
$p_1 \in {\mathfrak P}_1(I)$ and $p_2 \in {\mathfrak P}_2(I)$. 
Returning to (6.7)-(6.8) and applying Proposition 6.2 to $p_2$, we get 
$u_2(y) \leq I^{3/4}$, so
$$
u(x)^2 \leq I^{3/2} \int_{-\infty}^{+\infty} u_1(x-y)^2\, dy \leq I^{5/2}.
$$


\vskip10mm
\section{{\bf Bounds under moment assumptions}}
\setcounter{equation}{0}

\vskip2mm
\noindent
Another way to sharpen the bounds obtained in Section 2 for general densities 
with finite Fisher information is to invoke conditions on the absolute moments
$$
\beta_s = \beta_s(X) = \E\, |X|^s \qquad (s > 0 \ \ {\rm real}).
$$
By Proposition 2.1 and Cauchy's inequality, if the Fisher 
information is finite, 
\bee
\int_{-\infty}^{+\infty} |x|^s\, |p'(x)|\,dx
 & = &
\int_{\{p(x)>0\}} |x|^s p(x)^{1/2} \ \frac{|p'(x)|}{p(x)^{1/2}} \ dx \\
 & \leq & 
\bigg(\int_{\{p(x)>0\}} |x|^{2s} p(x)\,dx\bigg)^{1/2}\,
\bigg(\int_{\{p(x)>0\}} \frac{p'(x)^2}{p(x)}\,dx\bigg)^{1/2}.
\ene
Hence, we arrive at:

\vskip5mm
{\bf Proposition 7.1.} {\it If $X$ has an absolutely continuous density $p$,
then, for any $s>0$,
$$
\int_{-\infty}^{+\infty} |x|^s\, |p'(x)|\,dx \leq \sqrt{\beta_{2s} I(X)}.
$$
}

\vskip2mm
This bound holds irrespectively of the Fisher information or the $2s$-th 
absolute moment $\beta_{2s}$ being finite or not.

Below we describe several applications of this proposition.

First, let us note that, when $s \geq 1$, the function 
$u(x) = (1+|x|^s) p(x)$ is (locally) absolutely continuous and has a 
Radon-Nikodym derivative satisfying
$$
|u'(x)| \leq s |x|^{s-1}\,p(x) + (1+|x|^s)\, |p'(x)|. 
$$
Integrating this inequality and assuming that both $I(X)$ and $\beta_{2s}$ 
are finite, we see that $u$ is a function of bounded variation. Since $u$ 
is integrable as well, we have
$$
u(-\infty) = \lim_{x \rightarrow -\infty} u(x) = 0, \qquad
u(+\infty) = \lim_{x \rightarrow +\infty} u(x) = 0.
$$
Therefore, applying Propositions 2.2 and 7.1, we get
\bee
u(x) \, = \, \int_{-\infty}^x u'(y)\,dy 
 & \leq & 
\int_{-\infty}^{+\infty} |u'(y)|\,dy \\
 & \leq &
s\int_{-\infty}^{+\infty} |x|^{s-1}\,p(x)\,dx +
\int_{-\infty}^{+\infty} (1+|x|^s)\, |p'(x)|\,dx \\
 & \leq &
s \beta_{s-1} + \sqrt{I(X)} + \sqrt{\beta_{2s} I(X)}.
\ene
In addition, $u(x) \rightarrow 0$, as $x \rightarrow \infty$.
One can summarize.

\vskip5mm
{\bf Corollary 7.2.} {\it If $X$ has density $p$, then, given $s \geq 1$,
for any  $x \in \R$,
$$
p(x) \leq \frac{C}{1 + |x|^s}
$$
with a constant $C = s \beta_{s-1} + \sqrt{(1+\beta_{2s}) I(X)}$.
If this constant is finite, we also have
$$
\lim_{x \rightarrow \infty} \, (1 + |x|^s)\,p(x) = 0.
$$
}

\vskip2mm
In the resulting inequality no requirements on the density are needed.

Applying Proposition 7.1 and Corollary 7.2 (the last assertion) with $s=1$, 
we obtain the following sharpening of Corollary 2.3.

\vskip5mm
{\bf Corollary 7.3.} {\it If $X$ has finite second moment and finite Fisher 
information $I(X)$, then for its characteristic function 
$f(t) = \E\,e^{itX}$ we have
$$
|f'(t)| \leq \frac{C}{|t|}, \qquad t \in \R,
$$
with constant $C = 1 + \sqrt{\beta_2 I(X)}$.
}

\vskip5mm
Indeed, if $p$ is density of $X$ and $t \neq 0$, one may integrate by 
parts
\bee
f'(t)
 \ = \
\int_{-\infty}^{+\infty} e^{itx}\, (ix)\, p(x)\,dx 
 & = & 
\frac{1}{t} \int_{-\infty}^{+\infty} x p(x)\,d e^{itx} \\
 & = & 
- \frac{1}{t} \int_{-\infty}^{+\infty} (p(x) + x p'(x))\,e^{itx}\,dx,
\ene
which yields $|tf'(x)| \leq 1 + \sqrt{\beta_2 I(X)}$.

Under stronger moment assumptions, one can obtain better bounds in 
comparison with Corollary 7.2. For example, if for some 
$\lambda > 0$, the exponential moment
$$
\beta = \E\, e^{2\lambda |X|} = \int_{-\infty}^{+\infty} e^{2\lambda |x|}\,p(x)\,dx
$$
is finite, then by similar arguments, for any $x \in \R$, we have
$p(x) \leq C\,e^{-\lambda |x|}$
with some constant $C$ depending on $\lambda$, $\beta$ and $I(X)$.


\vskip10mm
\section{{\bf Fisher information in terms of the second derivative}}
\setcounter{equation}{0}

\vskip2mm
\noindent
It will be convenient to work with the formula for the Fisher information 
involving the second derivative of the density. We state it for convolutions 
of two densities with finite Fisher information.

\vskip5mm
{\bf Proposition 8.1.} {\it If a random variable $X$ has density 
$p \in {\mathfrak P}_2$, then
\be
I(X) = -\int_{-\infty}^{+\infty} p''(x)\,\log p(x)\,dx,
\en
provided that
\be
\int_{-\infty}^{+\infty} |p''(x)\,\log p(x)|\,dx < +\infty.
\en
The latter condition holds, if\, $\E\,|X|^s < +\infty$ for some $s > 2$.
}

\vskip5mm
Strictly speaking, the integration in (8.1)-(8.2) should be performed over 
the set $\{x:p(x)>0\}$. One may extend this integration to the whole real line
by using the convention $0 \log 0 = 0$. This is consistent with the property 
that $p''(x) = 0$, as soon as $p(x) = 0$ (according to Proposition 6.3).

\vskip2mm
{\bf Proof.} The assumption $p \in {\mathfrak P}_2$ ensures that $p$ has 
an absolutely continuous derivative $p'$ with Radon-Nikodym derivative
$p''$. By Proposition 6.2, $p'$ has bounded total variation, which
justifies the possibility of integration by parts.

More precisely, assuming that $p \in {\mathfrak P}_2$, let us decompose the 
open set $\{x: p(x) > 0\}$ into disjoint open intervals $(a_n,b_n)$, bounded 
or not. In particular, $p(a_n) = p(b_n) = 0$, and by the bound (6.5) of 
Proposition 6.2,
$$
|p'(x)\log p(x)| \leq \,I^{3/4} \sqrt{p(x)}\, |\log p(x)| \rightarrow 0, \quad
{\rm as} \ \ x \downarrow a_n,
$$
and similarly for $b_n$. Integrating by parts, we get for $a_n < T_1 < T_2 < b_n$,
\bee
\int_{T_1}^{T_2} \frac{p'(x)^2}{p(x)}\, dx
 & = &
\int_{T_1}^{T_2} p'(x)\, d\log p(x) \\
 & = &
p'(x) \log p(x) \bigg|_{x=T_1}^{T_2} - \int_{T_1}^{T_2} p''(x)\,\log p(x)\,dx.
\ene
Letting $T_1 \rightarrow a_n$ and $T_2 \rightarrow b_n$, we get
$$
\int_{a_n}^{b_n} \frac{p'(x)^2}{p(x)}\, dx = 
 - \int_{a_n}^{b_n} p''(x)\,\log p(x)\,dx,
$$
where the second integral is understood in the improper sense. It remains to 
perform summation over $n$ on the basis of (8.2), and then we obtain (8.1).

To verify the integrability condition (8.2), one may apply an integral 
bound of Proposition 6.3. Namely, using Cauchy's inequality,
for the integral in (8.2) we have
$$
\Big(\int_{\{p(x)>0\}} 
\frac{|p''(x)|}{\sqrt{p(x)}}\ \sqrt{p(x)}\, |\log p(x)|\,dx\Big)^2 \, \leq \,
I^2 \int_{-\infty}^{+\infty} p(x) \log^2 p(x)\,dx.
$$
If the moment $\beta_s = \E\, |X|^s$ is finite, Corollary 7.2 yields 
$$
p(x) \log^2 p(x) \, \leq \, C\, \frac{\log(e + |x|)}{1 + |x|^{s/2}}
$$
with constant $C$ depending on $I$ and $\beta_s$. The latter function is 
integrable in case $s>2$, so the integral in (8.2) is finite.
Proposition 8.1 is proved.

\vskip2mm
Of course, for smooth positive $p$, (8.1) remains valid without additional 
assumptions. However, then the integral should be understood in the improper 
sense (it exists and is finite, as long as $X$ has finite Fisher information).

In order to involve the standard moment assumption --  the finiteness 
of the second moment, we consider densities representable as 
convolutions of more than two densities with finite Fisher information.

\vskip5mm
{\bf Proposition 8.2.} {\it If a random variable $X$ has finite second moment
and density $p \in {\mathfrak P}_5$, then condition $(8.2)$ holds, and
$X$ has Fisher information given by $(8.1)$.
}

\vskip5mm
To show that $(8.2)$ is fulfilled, it suffices to prove the following pointwise 
bounds which are of independent interest.

\vskip5mm
{\bf Proposition 8.3.} {\it If\, $\E X^2 \leq 1$ and $X$ has density 
$p \in {\mathfrak P}_5(I)$, then with some absolute constant $C$, for all $x$,
\be
|p''(x)| \, \leq \, CI^3\, \frac{1}{1 + x^2}
\en
and
\be
|p''(x)\,\log p(x)| \, \leq \, CI^3\, \frac{\log(e + |x|)}{1 + x^2}.
\en
}

\vskip2mm
{\bf Proof.} The assumption $\E X^2 \leq 1$ implies $I \geq 1$ (by Cramer-Rao's
inequality). Also, the characteristic function $f(t) = \E\, e^{itX}$ 
is twice differentiable, and by Corollary 2.3, it satisfies
$$
|f(t)| \leq \frac{I^{5/2}}{|t|^5}.
$$
Hence, $p$ may be described as the inverse Fourier transform
$$
p(x) = \frac{1}{2\pi} \int_{-\infty}^{+\infty} e^{-itx} f(t)\,dt,
$$
and a similar representation is also valid for the second derivative,
\be
p''(x) = -\frac{1}{2\pi} \int_{-\infty}^{+\infty} e^{-itx}\, t^2 f(t)\,dt.
\en

Write $X = X_1 + \dots + X_5$ with independent summands such that $I(X_j) \leq I$ 
and assume (without loss of generality) that they have equal means. Then 
$\E X_j^2 \leq 1$, hence the characteristic functions $f_j(t)$ of $X_j$ have 
second derivatives $|f_j''(t)| \leq 1$. Moreover, by Corollaries 2.3 and 7.3, 
$$
|f_j(t)| \leq \frac{I^{1/2}}{|t|}, \qquad 
|f_j'(t)| \leq \frac{1 + I^{1/2}}{|t|}.
$$

Now, differentiation of the equality $f(t) = f_1(t) \dots f_5(t)$ leads to
$$
f'(t) = f_1'(t)\,f_2(t) \dots f_5(t) + \dots +
f_1(t) \dots f_4(t)\,f_5'(t),
$$
hence $|f'(t)| \leq \frac{5 I^2\, (1 + I^{1/2})}{|t|^5}$.
Differentiating once more, it should be clear that
$$
|f''(t)| \, \leq \, 
\frac{5 I^2}{t^4} + \frac{20\,I^{3/2} (1 + I^{1/2})^2}{|t|^5}.
$$
These estimates imply that 
$$
|(t^2 f(t))'| \leq \frac{CI^{5/2}}{|t|^3}, \qquad
|(t^2 f(t))''| \leq \frac{CI^{5/2}}{t^2} \qquad (|t| \geq 1)
$$
with some absolute constant $C$. As a consequence, one may differentiate
the equality (8.5) with $x \neq 0$ by parts to get
$$
p''(x) \, = \,
\frac{1}{2\pi\, (ix)^2} \int_{-\infty}^{+\infty} (t^2 f(t))''\,e^{-itx}\,dx.
$$
Hence, for all $x \in \R$,
\be
|p''(x)| \leq \frac{CI^{5/2}}{1 + x^2}
\en
with some absolute constant $C$.

Now, to derive the second pointwise bound, first we recall that 
$p(x) \leq I^{1/2}$. Hence,
\be
|\log p(x)| \leq \log(I^{1/2}) + \log \frac{I^{1/2}}{p(x)},
\en
where the last term is thus non-negative. Next, we partition the real 
line into the sets $A = \{x: p(x) \leq \frac{I^{1/2}}{2(1 + x^4)}\}$ and its 
complement $B$. On the set $A$, by Proposition 6.3,
$$
|p''(x)|\,\log \frac{I^{1/2}}{p(x)} \, \leq \, 
I^{5/4} \sqrt{p(x)}\,\log\frac{I^{1/2}}{p(x)}
 \, \leq \, C_1 I^{3/2}\, \frac{\log(e + |x|)}{1 + x^2},
$$
and similarly, by (8.6), on the set $B$ we have an analogous inequality
$$
|p''(x)|\,\log \frac{I^{1/2}}{p(x)} \, \leq \, |p''(x)|\,\log\big(2 (1 + x^4)\big)
 \, \leq \, C_2I^{5/2}\, \frac{\log(e + |x|)}{1 + x^2}.
$$
Thus, for all $x$, applying (8.7) and again (8.6),
\bee
|p''(x) \log p(x)| 
 & \leq & 
|p''(x)| \log(I^{1/2}) + |p''(x)| \log \frac{I^{1/2}}{p(x)} \\
 & \leq & 
C I^{5/2}\,(1 + \log I)\, \frac{\log(e + |x|)}{1 + x^2}.
\ene
Proposition 8.3 is proved.


\vskip10mm
\section{{\bf Normalized sums. Proof of Theorem 1.3}}
\setcounter{equation}{0}

\vskip2mm
\noindent
By the definition of classes ${\mathfrak P}_k$ ($k = 1,2,\dots$), 
the normalized sum
$$
Z_n = \frac{X_1 + \dots + X_n}{\sqrt{n}}
$$
of independent random variables $X_1,\dots,X_n$ with finite Fisher information
has density $p_n$ belonging to ${\mathfrak P}_k$, as long as $n \geq k$. 

Moreover, if all $I(X_j) \leq I$ for all $j$, then 
$p_n \in {\mathfrak P}_k(2kI)$. Indeed, one can partition the collection 
$X_1,\dots,X_n$ into $k$ groups and write $Z_n = U_1 + \dots + U_k$ with
$$
U_i = \frac{1}{\sqrt{n}}\,\sum_{j = i}^{m} X_{(i-1)m + j} \ \
(1 \leq i \leq k-1), \qquad
U_k = \frac{1}{\sqrt{n}}\,\sum_{j = (k-1)m + 1}^n X_j,
$$
where $m = [\frac{n}{k}]$. By Stam's inequality (6.1), for $1 \leq i \leq k-1$
$$
\frac{1}{I(U_i)} \, \geq \, \frac{1}{n}\,
\sum_{j = i}^{m} \frac{1}{I(X_{(i-1)m + j})} \, \geq \, \frac{m}{nI}
\, \geq \, \frac{1}{2kI},
$$
and similarly $\frac{1}{I(U_k)} \geq \frac{1}{2kI}$.

\vskip1mm
Therefore, the previous observations about densities from ${\mathfrak P}_k$ 
are applicable to $Z_n$ with sufficiently large $n$, as soon as the $X_j$ 
have finite Fisher information with a common bound on $I(X_j)$. 

A similar application of (6.1) also yields $I(Z_n) \leq 2 I(Z_{n_0})$.
Here, the factor $2$ may actually be removed, as a consequence of one
generalization of Stam's inequality obtained by Artstein, Ball, Barthe and Naor.
It is formulated below as a separate proposition
(although for our purposes the weaker inequality is sufficient).

\vskip5mm
{\bf Proposition 9.1} [A-B-B-N2]. {\it If $(X_n)_{n \geq 1}$ are independent 
and identically distributed, then 
$$
I(Z_n) \leq I(Z_{n_0}), \quad {for \ all} \ \ n \geq n_0.
$$
}

We are now ready to return to Theorem 1.3 and complete its proof.

\vskip5mm
{\bf Proof of Theorem 1.3.} Let $(X_n)_{n \geq 1}$ have finite second moment 
and a common characteristic function $f_1$. The characteristic function of
$Z_n$ is thus 
\be
f_n(t) = \E\, e^{itZ_n} = f_1\bigg(\frac{t}{\sqrt{n}}\bigg)^n.
\en

Clearly, $a) \Rightarrow b) \Leftrightarrow c)$. 

If $Z_n$ has density $p_n$ of bounded total variation, Proposition 4.1 
yields $I(Z_{3n}) = I(p_{3n}) \leq \frac{3}{2}\,\|p_n\|_{{\rm TV}}^2 < +\infty$.
Hence we obtain $c) \Rightarrow a)$, as well, and thus, the conditions $a)-c)$ 
are equivalent.

$a) \Rightarrow d)$. 
Assume that $I(Z_{n_0}) < +\infty$ for some fixed $n_0 \geq 1$. Applying 
Corollary 2.3 with $X = Z_{n_0}$, it follows that
$$
|f_{n_0}(t)| \leq \frac{1}{t}\,\sqrt{ n_0 I(Z_{n_0})}, \qquad t > 0.
$$
Hence, $|f_1(t)| \leq Ct^{-\ep}$ with constants $\ep = \frac{1}{n_0}$ and
$C = \big(n_0 I(Z_{n_0})\big)^{1/2n_0}$ which is $d)$.

$d) \Rightarrow e)$ is obvious.

$e) \Rightarrow c)$. Differentiating the formula (9.1) and using the 
integrability assumption (1.8) on $f_1$, we see that, for all 
$n \geq \nu + 2$, the characteristic function  $f_n$ and its first two 
derivatives are integrable with weight $|t|$. This implies in particular 
that $Z_n$ has a continuously differentiable density
\be
p_n(x) = \frac{1}{2\pi}\, \int_{-\infty}^{+\infty} e^{-itx} f_n(t)\,dt,
\en
which, by Proposition 5.1, has finite total variation
$$
\|p_n\|_{{\rm TV}} \, = \, \int_{-\infty}^{+\infty} |p_n'(x)|\,dx \, \leq \,
\frac{1}{2}\, 
\int_{-\infty}^{+\infty} \big(|tf_n''(t)| + 2\,|f_n'(t)| + |t f_n(t)|\big)\,dt.
$$
Thus, Theorem 1.3 is proved.

\vskip2mm
{\bf Remark 9.2.} 
If we assume in Theorem 1.3 finiteness of the first absolute moment of $X_1$ 
(rather than the finiteness of the second moment), the statement will remain 
valid, provided that the integrability condition $e)$ is replaced with 
a stronger condition like
\be
\int_{-\infty}^{+\infty} |f_1(t)|^\nu\,t^2\,dt < +\infty, \qquad 
{\rm for \ some} \ \ \nu > 0.
\en
In this case, it follows from (9.1) that, for all $n \geq \nu + 1$, the 
characteristic function  $f_n$ and its derivative are integrable with weight 
$t^2$. Therefore, according to Proposition 5.2, the normalized sum $Z_n$ has 
density $p_n$ with finite total variation
$$
\|p_n\|_{{\rm TV}} \, \leq \, \bigg(\int_{-\infty}^{+\infty} |t f_n(t)|^2\,dt
\int_{-\infty}^{+\infty} |(tf_n(t))'|^2\,dt\bigg)^{1/4}.
$$
As a result, we obtain the chain of implications
$(9.3) \Rightarrow b) \Rightarrow a) \Rightarrow d)$. The latter condition
ensures that $p_n$ admits the representation (9.2) and has a continuous 
derivative for sufficiently large $n$. That is, we obtain $c)$.


\vskip10mm
\section{{\bf Edgeworth-type expansions}}
\setcounter{equation}{0}

\vskip2mm
In the sequel, let $(X_n)_{n \geq 1}$ be independent identically distributed 
random variables with mean $\E X_1 = 0$ and variance $\Var(X_1) = 1$. 
Here we collect some auxiliary results about Edgeworth-type expansions 
for the distribution functions $F_n(x) = \P\{Z_n \leq x\}$ and the densities 
$p_n$ of the normalized sums $Z_n = (X_1 + \dots + X_n)/\sqrt{n}$.

If the absolute moment $\E\,|X_1|^s$ is finite for a given integer $s \geq 2$,
define
\be
\varphi_s(x) = \varphi(x) + \sum_{k=1}^{s-2} q_k(x)\,n^{-k/2}
\en
with the functions $q_k$ described in the introductory section, i.e.,
\be
q_k(x) \ = \, \varphi(x)\, \sum H_{k + 2j}(x) \,
\frac{1}{r_1!\dots r_k!}\, \bigg(\frac{\gamma_3}{3!}\bigg)^{r_1} \dots
\bigg(\frac{\gamma_{k+2}}{(k+2)!}\bigg)^{r_k}.
\en
Here, $H_k$ denotes the Chebyshev-Hermite polynomial of degree $k \geq 0$
with leading coefficient~1, and the summation runs over all non-negative 
solutions $(r_1,\dots,r_k)$ to the equation
$r_1 + 2 r_2 + \dots + k r_k = k$ with $j = r_1 + \dots + r_k$.

Put also
\be
\Phi_s(x) \, = \, \int_{-\infty}^x \varphi_s(y)\,dy
 \, = \, \Phi(x) + \sum_{k=1}^{s-2} Q_k(x)\,n^{-k/2}.
\en
Similarly to $q_k$, the functions $Q_k$ have an explicit description involving
the cumulants $\gamma_3,\dots,\gamma_{k+2}$ of $X_1$, namely,
$$
Q_k(x) \ = \, -\varphi(x) \sum H_{k + 2j-1}(x) \,
\frac{1}{r_1!\dots r_k!}\, \bigg(\frac{\gamma_3}{3!}\bigg)^{r_1} \dots
\bigg(\frac{\gamma_{k+2}}{(k+2)!}\bigg)^{r_k},
$$
where the summation is the same as in (10.2), cf. [B-RR] or [P]. 

The functions $\varphi_s$ and $\Phi_s$ are used
to approximate the density and distribution function of $Z_n$ with error 
of order smaller than $n^{-(s-2)/2}$. The following lemma is classical.

\vskip5mm
{\bf Lemma 10.1.} {\it Assume that\,
$\limsup_{|t| \rightarrow +\infty} |f_1(t)| < 1$. 
If $\E\,|X_1|^s < +\infty$ $(s \geq 3)$, then as $n \rightarrow \infty$,
uniformly over all $x$
\be
(1 +|x|^s) \big(F_n(x) - \Phi_{[s]}(x)\big) = o\big(n^{-(s-2)/2}\big).
\en
}

\vskip2mm
Let us emphasize that (10.4) remains valid for general real $s \geq 2$. Here, 
$\Phi_s$ should be replaced with $\Phi_{[s]}$. For the range $2 \leq s < 3$ 
the Cramer condition for the characteristic function is not used, and the 
result was obtained in [O-P]; the case $s \geq 3$
is treated in [P] (cf. Theorem 2, Ch.VI, p. 168).

We also need to describe the approximation of densities.
Recall that $Z_n$ have the characteristic functions
$$
f_n(t) = f_1\bigg(\frac{t}{\sqrt{n}}\bigg)^n,
$$
where $f_1$ stands for the characteristic function of $X_1$. If the Fisher 
information $I(Z_{n_0})$ is finite, then, by Corollary 2.3, 
$|f_{n_0}(t)| \leq \frac{c}{|t|}$ with some constant (namely, $c^2 = I(Z_{n_0})$). 
Hence, given $m \geq 1$, the characteristic functions of $Z_n$ admit 
a polynomial bound $|f_n(t)| \leq c_m\,|t|^{-m}$ for $n \geq m n_0$ and
with $c_m$ which does not depend on $t$. Thus, for all sufficiently large $n$, 
$Z_n$ have continuous bounded densities 
$$
p_n(x) = \frac{1}{2\pi}\, \int_{-\infty}^{+\infty} e^{-itx} f_n(t)\,dt,
$$
which have continuous derivatives 
\be
p_n^{(l)}(x) = \frac{1}{2\pi}\, \int_{-\infty}^{+\infty} 
(-it)^l\, e^{-itx} f_n(t)\,dt
\en
of any prescribed order.

\vskip5mm
{\bf Lemma 10.2.} {\it Assume $I(Z_{n_0}) < +\infty$, for some $n_0$,
and let $\E\, |X_1|^s < +\infty$ $(s \geq 2)$.
Fix $l = 0,1,\dots$ Then, for all sufficiently large $n$,
\be
(1 + |x|^s)\, |p_n^{(l)}(x) - \varphi_s^{(l)}(x)| \, \leq \, 
\psi_{l,n}(x)\,\frac{\ep_n}{n^{(s-2)/2}}, \qquad x \in \R,
\en
where $\ep_n \rightarrow 0$, as $n \rightarrow \infty$, and
\be
\sup_x \, |\psi_{l,n}(x)| \leq 1, \qquad 
\int_{-\infty}^{+\infty} \psi_{l,n}(x)^2\,dx \leq 1.
\en
}

\vskip5mm
In case $l=0$, this lemma with the first bound 
$\sup_x \, |\psi_{l,n}(x)| \leq 1$ is a well-known result, which does not 
need to require the finiteness of Fisher information, while using the 
assumption of the boundedness of $p_n$ for large $n$, only. We can refer to 
[P], p. 211 in case $s \geq 3$ and to [P], pp. 198-201 for the case $s=2$ 
when $\varphi_s = \varphi$. The result follows from the corresponding 
Edgeworth-type approximation of $f_n(t)$ by the Fourier transforms of 
$\varphi_s(x)$ on growing intervals such as $|t| < c_1 n^{1/6}$ in case 
$s \geq 3$. Repeating the arguments on pp. 211-212 of [P] and applying 
Plancherel's formula, one can easily obtain the second bound in (10.7), as well. 
In fact, the case $l \geq 1$ is similar, since the appearence of the additional 
factor $(-it)^l$ in (10.5) does not create any difficulty due to the polynomial 
decay at infinity of the characteristic functions $f_n$. 

For the proof of Theorem 1.1, the lemma will be used with the values $l = 0,1,2$, 
only.


\vskip10mm
\section{{\bf Behaviour of densities not far from the origin}}
\setcounter{equation}{0}

\vskip2mm
To study the asymptotic behavior of the Fisher information distance
$$
I(Z_n||Z) = 
\int_{-\infty}^{+\infty} \frac{(p_n'(x) + xp_n(x))^2}{p_n(x)}\ dx,
$$
we split the domain of integration into the interval $|x| \leq T_n$ and its
complement. Thus, define
$$
J_0 = \int_{|x| \leq T_n} \frac{(p_n'(x) + xp_n(x))^2}{p_n(x)}\ dx
$$
and similarly $J_1$ for the region $|x| > T_n$.
If $T_n$ is not too large, the first integral can be treated with the help of 
Lemma 10.2. Namely, we take
\be
T_n = \sqrt{(s-2)\log n + s \log \log n + \rho_n} \qquad (s>2),
\en
where $\rho_n \rightarrow +\infty$ is a sufficiently slowly growing sequence 
whose growth is restricted by the decay of the sequence $\ep_n$ in (10.6). 
In other words, $[-T_n,T_n]$ represents an asymptotically largest interval, 
where we can guarantee that the densities $p_n$ of $Z_n$ are separated from 
zero, and moreover, 
$\sup_{|x| \leq T_n} |\frac{p_n(x)}{\varphi(x)} - 1| \rightarrow 0$.
To cover the case $s=2$, one may put $T_n = \sqrt{\rho_n}$, where
$T_n \rightarrow +\infty$ is a sufficiently slowly growing sequence.
With this choice of $T_n$, an estimation of the integral $J_1$ can be
performed via moderate inequalities.

In this section we focus on $J_0$ and provide an asymptotic expansion for 
it with a remainder term which turns out to be slightly better in comparison 
with the resulting expansion (1.3) of Theorem~1.1.

\vskip5mm
{\bf Lemma 11.1.} {\it Let $s \geq 3$ be an integer. If $I(Z_{n_0}) < +\infty$, 
for some $n_0$, then
$$
J_0 = 
\frac{c_1}{n} + \frac{c_2}{n^{2}} + \dots + \frac{c_{[(s-2)/2]}}{n^{[(s-2)/2]}}
+ o\Big(\frac 1{n^{(s-2)/2}\,(\log n)^{(s-1)/2}}\Big),
$$
where the coefficients $c_j$ are defined in $(1.4)$.
}

\vskip5mm
{\bf Proof}. Let us adopt the convention to write $\delta_n$ for any sequence
of functions satisfying $|\delta_n(x)| \leq \ep_n n^{-(s-2)/2}$ with
$\ep_n \rightarrow 0$, as $n \rightarrow \infty$, at least on the intervals
$|x| \leq T_n$. For example, the statement of Lemma 10.2 with $l=0$
may be written as
\be
p_n(x) = (1 + u_s(x))\varphi(x) + \frac{\delta_n}{1+|x|^s},
\en
where 
$$
u_s(x) \, = \, \frac{\varphi_s(x)-\varphi(x)}{\varphi(x)}
\, = \, \sum_{k=1}^{s-2}\, \frac{q_k(x)}{\varphi(x)} \ \frac{1}{n^{k/2}}.
$$
Combining the lemma with $l=0$ and $l=1$, we obtain another representation
\be
p_n'(x) + xp_n(x) = w_s(x) + \frac{\delta_n}{1+|x|^{s-1}},
\en
where 
$$
w_s(x) \, = \, \sum_{k=1}^{s-2}\, \frac{q'_k(x)+xq_k(x)}{n^{k/2}}.
$$

Note that the functions $u_s$ and $w_s$ depend on $n$ as parameter and 
are getting small for growing $n$. More precisely, it follows from the 
definition of $q_k$ that, for all $x \in \R$,
\be
\frac{|w_s(x)|}{\varphi(x)} \, \leq \, C_s\frac{1+|x|^{3(s-1)}}{\sqrt n}
\qquad \text{and} \qquad
|u_s(x)| \, \leq \, C_s\frac{1+|x|^{3(s-2)}}{\sqrt n} 
\en
with some constants depending on $s$ and the cumulants of $X_1$, only. 
In particular, for $|x| \leq T_n$ and any prescribed $0 < \ep < \frac{1}{2}$,
\be
\frac{|w_s(x)|}{\varphi(x)} \, < \, \frac{1}{n^{\frac{1}{2} - \ep}} 
\qquad \text{and} \qquad
|u_s(x)| \, < \, \frac{1}{4}
\en
with sufficiently large $n$. In addition, with a properly chosen sequence 
$\rho_n$, we have
\be
\frac{\delta_n}{T_n^s\,\varphi(T_n)} \, < \, \frac{1}{4}.
\en
Hence, by Lemma 10.2, $|\frac{p_n(x)}{\varphi(x)} - 1| < \frac{1}{2}$ 
on the interval $|x| \leq T_n$.

Now, for $|x| \leq T_n$
$$
\big(1+u_s(x)\big)^{-1} - 
\Big(1 + u_s(x) + \frac{\delta_n}{(1+|x|^s)\varphi(x)}\Big)^{-1}
 = \frac{\delta_n}{(1+|x|^s)\varphi(x)},
$$
and we obtain from (11.2)
$$
\frac{1}{p_n(x)} \, = \, \frac{1}{(1+u_s(x))\varphi(x)} + 
\frac{\delta_n}{(1+|x|^s)\varphi(x)^2}. 
$$
Combining this with (11.3) and using (11.5), we will be lead to
$$
\frac{(p_n'(x)+xp_n(x))^2}{p_n(x)} = \frac{w_s(x)^2}
{(1+u_s(x))\varphi(x)}+\sum_{j=1}^5 r_{nj}(x), \qquad |x| \leq T_n,
$$
where
\bee
r_{n1} & = &
\frac{w_s(x)}{(1+|x|^{s-1})\varphi(x)} \ \delta_n, \qquad \ \ 
r_{n2} \ = \
\frac {w_s(x)^2}{(1+|x|^s)\varphi(x)^2} \ \delta_n, \\
r_{n3} & = &
\frac {w_s(x)}{(1+|x|^{{2s-1}})\varphi(x)^2} \ \delta_n^2, \qquad
r_{n4} \ = \
\frac{1}{(1+|x|^{2s-2})\varphi(x)} \ \delta_n^2, \\
r_{n5} & = &
\frac 1{(1+|x|^{3s-2})\varphi(x)^2} \ \delta_n^3.
\ene
Here, according to the left inequality in (11.5), the remainder terms
$r_{n1}(x)$ and $r_{n2}(x)$ are uniformly bounded on $[-T_n,T_n]$ 
by $|\delta_n|\, n^{-1/3}$. A similar bound also holds for $r_{n3}(x)$,
by taking into account (11.6). In addition, integrating by parts,
for large $n$ and with some constants (independent of $n$), we have
\bee
\int_{|x|\le T_n}|r_{n4}(x)|\,dx 
 & \leq &
\frac{C\ep_n}{n^{s-2}} \int_1^{T_n} \frac{1}{x^{2s-2}}\,e^{x^2/2}\,dx  \\
 & \leq & 
\frac{C' \ep_n}{n^{s-2}} \ \frac{1}{T_n^{2s-1}} \, e^{T_n^2/2}
 \ = \
o\Big(\frac {1}{T_n^{s - 1}\, n^{(s-2)/2}}\Big).
\ene
With a similar argument, the same $o$-relation also holds for the integral 
of $|r_{n5}(x)|$.

Thus,
\be
\int_{|x|\le T_n} \frac{(p_n'+xp_n)^2}{p_n} \ dx =
\int_{|x|\le T_n} \frac{w_s^2}{(1+u_s)\varphi} \ dx + 
o\Big(\frac {1}{T_n^{s-1}n^{(s-2)/2}}\Big). 
\en

Now, by Taylor's expansion around zero, in the interval $|u|\le \frac{1}{4}$
we have
$$
\frac{1}{1+u} \, = \, \sum_{k=0}^{s-4}\, (-1)^ku^k + \theta u^{s-3}, \qquad 
|\theta| < 2
$$
(there are no terms in the sum for $s=3$). Hence, with some $-2 < \theta_n < 2$
$$
\int_{|x|\le T_n}\frac{w_s^2}{(1+u_s)\varphi}\,dx \, = \,
\sum_{k=0}^{s-4}\, (-1)^k \int_{|x|\le T_n} w_s^2 u_s^k\,\frac{dx}{\varphi}
+ \theta_n \int_{|x|\le T_n} w_s^2 u_s^{s-3}\,\frac{dx}{\varphi}. 
$$
At the expense of a small error, these integrals may be extended to the whole
real line. Indeed, for large enough $n$, by (11.4), we have, for 
$k=0,1,\dots,s-4$ with some common constant $C_s$
$$
\int_{|x|>T_n} w_s^2\, |u_s|^k\,\frac{dx}{\varphi}
 \, \leq \,
\frac{C_s}{n^{(k+2)/2}}\int_{|x|>T_n}(1+|x|^{(3k+6)(s-1)})\, \varphi\,dx
 \, = \,
o\Big(\frac{1}{n^{(s-1)/2}}\Big).
$$
Moreover,
$$
\int_{-\infty}^{+\infty} w_s^2\, |u_s|^{s-3}\,\frac{dx}{\varphi} 
 \, = \, O\Big(\frac{1}{n^{(s-1)/2}}\Big).
$$
Therefore,
$$
\int_{|x|\le T_n} \frac{w_s^2}{(1+u_s)\varphi}\,dx \, = \, 
\sum_{k=0}^{s-4} \, (-1)^k
\int_{-\infty}^{+\infty}w_s^2 u_s^k\,\frac{dx}{\varphi}
+ O\Big(\frac{1}{n^{(s-1)/2}}\Big).
$$
Inserting this in (11.7), we thus arrive at
\be
J_0 \, = \, \sum_{k=0}^{s-4} \, (-1)^k \int_{-\infty}^{+\infty}
w_s^2 u_s^k\,\frac{dx}{\varphi} + o\Big(\frac{1}{T_n^{s-1}n^{(s-2)/2}}\Big).
\en

In the next step, we develop this representation by expressing 
$u_s$ and $w_s$ in terms of $q_k$ while expanding the sum in (11.8) 
in powers of $1/\sqrt{n}$ as
$$
\sum_{j=2}^{s-2} \, \frac{a_j}{n^{j/2}} + O\Big(\frac{1}{n^{(s-1)/2}}\Big).
$$

More precisely, here the coefficients are given by
\be
a_j \, = \, \sum_{k=2}^{j}\, (-1)^k \int_{-\infty}^{+\infty}
(q_{r_1}' + xq_{r_1})\, (q_{r_2}' + xq_{r_2})\,
q_{r_3} \dots q_{r_k}\ \frac{dx}{\varphi^{k - 1}} 
\en
with summation over all positive solutions $(r_1,\dots,r_k)$ to
$r_1 + \dots + r_k = j$. Moreover, when $j$ are odd, the above integrals 
are vanishing. Indeed, differentiating the equality (10.2) 
which defines the functions $q_k$ and using the property 
$H_n'(x) = n H_{n-1}(x)$ $(n \geq 1)$, we obtain a similar equality
\be
q_k'(x) + xq_k(x) \ = \, \varphi(x)\, \sum (k + 2l)\,H_{k + 2l - 1}(x) \,
\frac{1}{r_1!\dots r_k!}\, \bigg(\frac{\gamma_3}{3!}\bigg)^{r_1} \dots
\bigg(\frac{\gamma_{k+2}}{(k+2)!}\bigg)^{r_k}
\en
with summation over all non-negative solutions $(r_1,\dots,r_k)$ to
$r_1 + 2 r_2 + \dots + k r_k = k$, and where $l = r_1 + \dots + r_k$.
Hence, the integrand in (11.9) represents a linear combination of the
functions of the form
$$
H_{r_1 + 2l_1 - 1}\, H_{r_2 + 2l_2 - 1}\, H_{r_3 + 2l_3} \dots H_{r_k + 2l_k}\,
\varphi.
$$
Note that here the sum of indices is ${{\rm mod}\, 2}$ the same as $j$.
We can now apply the following property of the Chebyshev-Hermite polynomials 
(see Szeg\"o 1967). If the sum of indices $d_1,\dots, d_k$ is odd, 
then necessarily
$$
\int_{-\infty}^{\infty} H_{d_1}(x) \dots H_{d_k}(x) \, \varphi(x)\,dx = 0.
$$

Hence, $a_j = 0$, whenever $j$ is odd, and putting $c_j = a_{2j}$,
we arrive at the assertion of the lemma.

\vskip5mm
{\bf Remark.} In formula (11.9) with $c_j = a_{2j}$ we perform summation
over all integers $r_l \geq 1$ such that $r_1 + \dots + r_k = 2j$.
Hence, all $r_l \leq 2j - 1$, and thus the functions $q_{r_l}$ are determined
by the cumulants up to order $2j+1$. Hence, $c_j$ represents a polynomial in
$\gamma_3,\dots,\gamma_{2j+1}$.


\vskip10mm
\section{{\bf Moderate deviations}}
\setcounter{equation}{0}

\vskip2mm
We now consider the second integral
$$
J_1 = 
\int_{|x|>T_n} \frac{(p_n'(x)+xp_n(x))^2}{p_n(x)}\,dx
$$
participating in the Fisher information distance $I(Z_n||Z)$.

\vskip5mm
{\bf Lemma 12.1.} {\it Let $s \geq 3$ be an integer. If $I(Z_{n_0}) < +\infty$, 
for some $n_0$, then 
$$
J_1 = o\Big(\frac 1{n^{(s-2)/2}(\log n)^{(s-3)/2}}\Big).
$$
}

\vskip2mm
{\bf Proof.} Write
\be
J_1 \, \le \, 2J_{1,1} + 2J_{1,2} \, = \, 
2\int_{|x|>T_n} \frac{p_n'(x)^2}{p_n(x)}\,dx + 2\int_{|x|>T_n} x^2p_n(x)\,dx.
\en
Using Lemma 10.1, we conclude that, for $s=3,\dots$,
\be
J_{1,2} = o\Big(\frac 1{(n\log n)^{(s-2)/2}}\Big). 
\en
Indeed, integrating by parts we have
$$
\int_{T_n}^{+\infty} x^2 p_{n}(x)\,dx \, = \, T_n^2\,(1-F_n(T_n)) +
2\int_{T_n}^{+\infty} x(1-F_n(x)) \,dx.
$$
Recalling the definition (10.3) of the approximating functions 
$\Phi_s$ and applying an elementary inequality 
$1-\Phi(x) < \frac{1}{x}\,\varphi(x)$ ($x>0$), we obtain from (10.4)
\bee
T_n^2\, (1-F_n(T_n))
 & = &
T_n^2\, (1-\Phi_s(T_n)) + T_n^2\, (\Phi_s(T_n)-F_n(T_n)) \\
 & \leq & 
T_n \varphi(T_n) + C\,\varphi(T_n)\, \sum_{k=1}^{s-2}\, T_n^{3k} n^{-k/2}
+ o\Big(\frac{1}{T_n^{s-2}\, n^{(s-2)/2}}\Big) \\
 & = &
o\Big(\frac{1}{(n \log n)^{(s-2)/2}}\Big)
\ene
with some constant $C$. In addition,
\bee
\int_{T_n}^{+\infty} x(1-F_n(x)) \,dx
 & \leq & 
1-\Phi(T_n) +
C \sum_{k=1}^{s-2}\frac 1{n^{k/2}}\int_{T_n}^{+\infty} x^{3k}\varphi(x)\,dx \\
 & &
+ \, o\Big(\frac{1}{T_{n}^{s-2}n^{(s-2)/2}}\Big) \ = \
o\Big(\frac{1}{(n \log n)^{(s-2)/2}}\Big).  
\ene
With similar estimates for the half-axis $x<-T_n$, we arrive at the relation (12.2).

Let us now estimate $J_{1,1}$. Denote by $J_{1,1}^+$ the part of this integral
corresponding to the interval $x > T_n$. By Propositions 6.2, 6.4 and 8.3, 
for sufficiently large $n$ one may integrate by parts to justify the formula
\be
J_{1,1}^+ = - p_n'(T_n) \log p_n(T_n) -
\int_{T_n}^{+\infty} p_n''(x)\log p_n(x)\,dx. 
\en
Since $p_n(x) \leq C\sqrt{I(Z_{n_0})}$ for all $x$ (Propositions 2.2 and 9.1)
and since $p_n(T_n) \geq \frac{1}{2}\,\varphi(T_n)$, we see that 
for all sufficiently large $n$, $|\log p_n(T_n)|\le c T_n^2$ with some 
constants $C$ and $c$. Therefore, by Lemma 10.2 for the derivative of 
the density $p_n$, we get
\begin{eqnarray}
|p_n'(T_n) \log p_n(T_n)|
 & \leq & 
c T_n^2\, |p_n'(T_n)| \nonumber \\ 
 & \leq & 
cT_n^2\, |\varphi'(T_n)| + o\Big(\frac 1{T_n^{s-2}\, n^{(s-2)/2}}\Big) \nonumber \\
 & = &
o\Big(\frac 1{T_n^{s-3}\, n^{(s-2)/2}}\Big).
\end{eqnarray}
A similar relation holds at the point $-T_n$, as well.

It remains to evaluate the integral in (12.3). First we integrate over the set
$A = \{x > T_n: p_n(x)\le\varphi(x)^4\}$. By the upper bound of Proposition 6.4 
and applying Proposition 9.1 once more, we have, for all $x$ and all 
sufficiently large $n$, with some constant $C$
$$
|p_n''(x)| \, \leq \, I(p_n)^{5/4} \sqrt{p_n(x)} \, \leq \,
C I(Z_{n_0})^{5/4} \sqrt{p_n(x)}.
$$
Hence, with some constants $c,c'$
\bee
\int_A |p_n''(x)\log p_n(x)|\,dx
 & \leq & 
c \int_A\sqrt{p_n(x)}\, |\log p_n(x)|\,dx \\
 & \leq & 
c'\int_{T_n}^{+\infty} x^2\varphi(x)^2\,dx \, = \, o\Big(\frac{1}{n^{s-2}}\Big).  
\ene

On the other hand, for the complementary set $B = (T_n,+\infty) \setminus A$,
we have
\begin{equation}
\int_B |p_n''(x)\log p_n(x)|\,dx \, \le \, c\int_B x^2\, |p_n''(x)|\,dx. 
\end{equation}
We now apply Lemma 10.2 to approximate the second derivative. It yields
$$
\int_{T_n}^{+\infty} x^2\, |p_n''(x)|\,dx \, \leq \,
\int_{T_n}^{+\infty} x^2\, |\varphi_s''(x)|\,dx + \int_{T_n}^{+\infty} 
\frac{|\psi_{2,n}(x)|}{1+|x|^{s-2}}\,dx \cdot o\Big(\frac{1}{n^{(s-2)/2}}\Big).
$$
Here, the first integral on the right-hand side is bounded by
$$
\int_{T_n}^{+\infty} x^2\, |\varphi_s''(x)-\varphi''(x)|\,dx + 
\int_{T_n}^{+\infty} x^2\, |x^2-1|\,\varphi(x)\,dx = 
o\Big(\frac 1{T_n^{s-3}n^{(s-2)/2}}\Big).
$$
To estimate the second integral, we use Cauchy's inequality, which gives
$$
\int_{T_n}^{+\infty} \frac{1}{1 + |x|^{s-2}} \, |\psi_{2,n}(x)|\,dx
 \ \leq \ \frac{1}{T_n^{s-5/2}}  \
\bigg(\int_{-\infty}^{+\infty} \psi_{2,n}(x)^2\,dx\bigg)^{1/2}
 \ \leq \ \frac{1}{T_n^{s-5/2}}.
$$
Therefore, returning to (12.5), we get
$$
\int_B |p_n''(x)\log p_n(x)|\,dx \, = \, 
o\Big(\frac{1}{n^{(s-2)/2}\,(\log n)^{(s-3)/2}}\Big).
$$
Together with the bound for the integral over the set $A$, we thus have
$$
J_{1,1}^+ = o\Big(\frac{1}{n^{(s-2)/2}\,(\log n)^{(s-3)/2}}\Big).
$$

The part of the integral $J_{1,1}$ taken over the axis $x < -T_n$
admits a similar bound, hence the lemma is proved.

The statement of Theorem 1.1 in case $s \ge 3$ thus follows from 
Lemmas 11.1 and 12.1.


\vskip10mm
\section{{\bf Theorem 1.1 in the case $s=2$ and Corollary 1.2}}
\setcounter{equation}{0}

\vskip2mm
In the most general case $s=2$ the proof of Theorem 1.1 does no need
Edgeworth-type expansions. With tools developed in the previous sections 
the argument is straightforward and may be viewed as an alternative 
approach to Barron-Johnson's theorem.

To give more details, recall that once the Fisher information $I(Z_{n_0})$ 
is finite, the normalized sums $Z_n$ with $n \geq n_0 + 1$ have uniformly 
bounded densities $p_n$ with bounded continuous derivatives $p_n'$ 
(Proposition 6.2). Moreover, we have a well-known local limit theorem 
for densities; we described one of its variants in Lemma 10.2. 
In particular,
\begin{eqnarray}
\sup_x\ (1+x^2)\,|p_n(x)-\varphi(x)| & = & o(1), \\
\sup_x\ (1+x^2)\,|p_n'(x)-\varphi'(x)| & = & o(1),
\end{eqnarray}
as $n \to\infty$, where the convergence of the derivatives relies upon 
the finiteness of the Fisher information.

Splitting the integration in
$$
I(Z_n||Z) = 
\int_{-\infty}^{+\infty} \frac{(p_n'(x) + xp_n(x))^2}{p_n(x)}\ dx
$$
into the two regions, we have therefore, for every fixed $T>1$,
\be
J_0 \, = \,
\int_{|x|\le T}\frac{(p_n'(x)+xp_n(x))^2}{p_n(x)}\,dx \, = \,o(1),
\quad n\to\infty.
\en
On the other hand, write as we did before
\bee
J_1
   & = &
\int_{|x|>T}\frac{(p_n'(x)+xp_n(x))^2}{p_n(x)}\,dx  \, \leq \,
2 J_{1,1} + 2J_{1,2} \\ 
 & = & 
2\int_{|x|>T}\frac{p_n'(x)^2}{p_n(x)}\,dx + 2\int_{|x|>T} x^2 p_n(x)\,dx.
\ene
As we saw in (12.3),
$$
J_{1,1} = -p_n'(T)\log p_n(T) + p_n'(-T)\log p_n(-T)-
\int_{|x|>T} p_n''(x)\log p_n(x)\,dx. 
$$
By (13.1)-(13.2), $|p_n'(\pm T)\log p_n(\pm T)|\le 2T^3e^{-T^2/2}$
for all sufficiently large $n \geq n_T$. By Proposition 8.3, with some 
constant $c$, for all $x$,
$$
|p_n''(x) \log p_n(x)| \, \leq \, c\, \frac{\log(e+|x|)}{1+x^2},
$$
implying
$$
\int_{|x|>T} |p_n''(x)\log p_n(x)|\,dx \, \leq \, c' T^{-1/2} 
$$
with some other constant $c'$. In addition, by (13.1),
\bee
\int_{|x|>T} x^2 p_n(x)\,dx 
 & = &
\int_{|x|>T} x^2 (p_n(x) - \varphi(x))\,dx  + \int_{|x|>T} x^2 \varphi(x)\,dx \\
 & \hskip-20mm = & \hskip-10mm
-\int_{|x| \leq T} x^2 (p_n(x) - \varphi(x))\,dx  + \int_{|x|>T} x^2 \varphi(x)\,dx \\
 & \hskip-20mm \leq & \hskip-10mm
\int_{|x| \leq T} x^2\, |p_n(x) - \varphi(x)|\,dx  + \int_{|x|>T} x^2 \varphi(x)\,dx
\ \leq \ 2T^3\, o(1) + 4T\varphi(T).
\ene
Hence, given $\ep>0$, one can choose $T$ such that $J_1 < \ep$, for all 
$n$ large enough. This means that $J_1 = o(1)$, and recalling (13.3), 
we get $I(Z_n||Z) = o(1)$.

Let us now return to the case $s \geq 3$. 

\vskip2mm
{\bf Proof of Corollary 1.2.}
According to the expansion (11.8) which appeared in the proof of Lemma 11.1, 
Theorem 1.1 may equivalently be formulated as
\be
I(Z_n||Z) \, = \, \sum_{l=0}^{s-4} \, (-1)^l \int_{-\infty}^{+\infty}
w_s(x)^2 u_s(x)^l\,\frac{dx}{\varphi(x)} + 
o\Big(\frac{1}{n^{(s-2)/2} \, (\log n)^{(s-3)/2}}\Big),
\en
where as before
$$
w_s(x) \, = \, \sum_{j=1}^{s-2}\, (q'_j(x)+xq_j(x))\,n^{-j/2}, \qquad
u_s(x) \, = \, \sum_{j=1}^{s-2}\, \frac{q_j(x)}{\varphi(x)} \, n^{-j/2}.
$$

This representation for the Fisher information distance is more convenient 
for applications such as Corollary 1.2 in comparison with (1.3).
Assume that $s \geq 4$ and $\gamma_3 = \dots = \gamma_{k-1} = 0$ for
a given integer $3 \leq k \leq s$ (with no restriction when $k = 3$).
Then, by the definition (10.2), $q_1 = \dots = q_{k-3} = 0$, so
\be
w_s(x) \, = \, \sum_{j=k-2}^{s-2}\, (q'_j(x)+xq_j(x))\,n^{-j/2}, \qquad
u_s(x) \, = \, \sum_{j=k-2}^{s-2}\, \frac{q_j(x)}{\varphi(x)} \, n^{-j/2}.
\en
Hence, in order to isolate the leading term in (1.3) with the smallest power 
of $1/n$, one should take $l = 0$ in (13.4) and $j = k-2$ in the first sum of 
(13.5). This gives
\bee
I(Z_n||Z) & = & n^{-(k-2)}
\int_{-\infty}^{+\infty} 
\big(q'_{k-2}(x) + x q_{k-2}(x)\big)^2\,\frac{dx}{\varphi(x)} \\
 & & \hskip20mm + \
O\big(n^{-(k-1)}\big) + o\Big(\frac{1}{n^{(s-2)/2} \, (\log n)^{(s-3)/2}}\Big).
\ene
Now, again according to (10.2), or as found in (11.10),
$$
q'_{k-2}(x) + x q_{k-2}(x) = \frac{\gamma_k}{(k-1)!}\,H_{k-1}(x)\, \varphi(x).
$$
Therefore, the sum in (1.3) will contain powers of $1/n$ starting from
$1/n^{k-2}$ with leading coefficient
$$
c_{k-2} = \frac{\gamma_k^2}{(k-1)!^{\,2}}\,\int_{-\infty}^{+\infty} 
H_{k-1}(x)^2\, \varphi(x)\,dx = \frac{\gamma_k^2}{(k-1)!}.
$$
Thus, $c_1 = \dots = c_{k-3} = 0$ and we get
$$
I(Z_n||Z) \, = \, 
\frac{\gamma_k^2}{(k-1)!}\, \frac{1}{n^{k-2}} + O\big(n^{-(k-1)}\big) + o\Big(\frac{1}{n^{(s-2)/2} \, (\log n)^{(s-3)/2}}\Big).
$$


\vskip10mm
\section{{\bf Extensions to non-integer $s$. Lower bounds}}
\setcounter{equation}{0}

\vskip2mm
If $s \geq 2$ is not necessary integer, put $m=[s]$ (integer part).
Theorem 1.1 admits the following generalization. As before, let the normalized
sums
$$
Z_n = \frac{X_1 + \dots + X_n}{\sqrt{n}}
$$
be defined for independent identically distributed random variables with mean
$\E X_1=0$ and variance $\Var(X_1)=1$.

\vskip5mm
{\bf Theorem 14.1.} {\it If $I(Z_{n_0}) < +\infty$ for some $n_0$, and 
$\E\, |X_1|^s < +\infty$ $(s>2)$, then
\be
I(Z_n||Z) = \frac{c_1}{n} + \frac{c_2}{n^{2}} + \dots +
\frac{c_{[(s-2)/2]}}{n^{[(s-2)/2]}} +
o\Big(\frac 1{n^{(s-2)/2}\,(\log n)^{(s-3)/2}}\Big), 
\en
where the coefficients $c_j$ are the same as in $(1.4)$.
}

\vskip5mm
The proof is based on a certain extension and refinement
of the local limit theorem described in Lemma 10.2.

\vskip5mm
{\bf Lemma 14.2.} {\it Assume that $I(Z_{n_0}) < +\infty$ for some $n_0$,
and $\E\, |X_1|^s < +\infty$ $(s \geq 2)$. Fix $l = 0,1,\dots$ 
Then for all $n$ large enough, $Z_n$ have densities $p_n$ of class $C^l$ 
satisfying, as $n \to \infty$,
\be
(1+|x|^m)\, (p_n^{(l)}(x)-\varphi_m^{(l)}(x)) \, = \,
\psi_{l,n}(x)\,o(n^{-(s-2)/2}) 
\en
uniformly for all $x$, with\, $\sup_x\, |\psi_{l,n}(x)|\le 1$ and 
$\int_{-\infty}^{+\infty} \psi_{l,n}(x)^2\, dx \leq 1$. Moreover,
uniformly for all $x$, 
\begin{eqnarray}
\hskip-10mm (1+|x|^s)\,(p_n^{(l)}(x)-\varphi_m^{(l)}(x))
 & = &
\psi_{l,n,1}(x)\,o(n^{-(s-2)/2}) \nonumber\\
 & &
\hskip-20mm + \ (1+|x|^{s-m})\,\psi_{l,n,2}(x)\,
\big(O(n^{-(m-1)/2}) + o(n^{-(s-2)})\big),
\end{eqnarray}
where\, $\sup_x\, |\psi_{l,n,j}(x)| \leq 1$ and 
$\int_{-\infty}^{+\infty} \psi_{l,n,j}(x)^2\, dx \leq 1$ $(j=1,2)$.
}

\vskip5mm
Here we use the approximating functions 
$\varphi_m = \varphi + \sum_{k=1}^{m-2} q_k\, n^{-k/2}$ as before.

When $l = 0$ and in a simpler form, namely, with $\psi_{l,s,j}(x,n) = 1$, 
this result has recently been obtained in [B-C-G1]. In this case,
the finiteness of the Fisher information may be relaxed to the boundedness
of the densities. The more general case involving derivatives can be carried
out by a similar analysis as that developed in [B-C-G1], so we omit details.

If $s=m$ is integer, the Edgeworth-type expansions (14.2) and (14.3) coincide,
and we are reduced to the statement of Lemma 10.2. However, if $s>m$, (14.3) 
gives an improvement over (14.2) on relatively large intervals such as 
$|x| \leq T_n$ considered in Theorem 1.1 and defined in (11.1). 

\vskip2mm
{\bf Proof of Theorem 14.1.} With a few modifications one can argue in the 
same way as we did in the proof of Theorem 1.1. 
First, in case $l=0$ (14.3) yields, uniformly in $|x| \leq T_n$
$$
p_n(x) \, = \, \varphi_m(x) + \frac{1}{1+|x|^s}\, o\big(n^{-(s-2)/2}\big),
$$
which being combined with a similar relation for the derivative $(l=1)$ yields
$$
p_n'(x) + x p_n(x) \, = \, 
w_m(x) + \frac{1}{1+|x|^{s-1}}\, o\big(n^{-(s-2)/2}\big),
$$
where $w_m(x) = \sum_{k=1}^{m-2}\, (q'_k(x)+xq_k(x))\,n^{-k/2}$.
These two relations thus extend (11.2) and (11.3) which were only needed
in the proof of Lemma 11.1. Repeating the same arguments using
the functions $u_m(x) = \frac{\varphi_m(x)-\varphi(x)}{\varphi(x)}$,
we can extend the expansion of Lemma 11.1 with the same remainder term
to general values $s > 2$.

In order to prove Lemma 12.1 with real $s>2$, let us return to (12.1).
The fact that the relation (12.2) extends to non-integer $s$ follows from
the extended variant of Lemma 10.1, which was already mentioned before.
Thus our main concern has to be the integral $J_{1,1}$ which is responsible 
for the most essential contribution in the resulting remainder term.
Thus, consider the part of this integral on the positive half-axis
\be
J_{1,1}^+ =
\int_{T_n}^{+\infty} \frac{p_n'(x)^2}{p_n(x)}\,dx = 
- p_n'(T_n) \log p_n(T_n) - \int_{T_n}^{+\infty} p_n''(x)\log p_n(x)\,dx.
\en

Applying (14.3) at $x=T_n$, we obtain (12.4) for real $s>2$, that is,
$$
|p_n'(T_n) \log p_n(T_n)| = o\Big(\frac{1}{n^{(s-2)/2}\, (\log n)^{s-3}}\Big).
$$

To prove (14.1), it remains to estimate the last integral in (14.4) which has 
to be treated with an extra care. The argument uses both (14.2) and (14.3)
which are applied on different parts of the half-axis $x>T_n$.
For the set $A = \{x \geq T_n: p_n(x) \leq \varphi(x)^4\}$ we have already 
obtained a general relation
$$
\int_A |p_n''(x)\log p_n(x)|\,dx \, = \, o\Big(\frac{1}{n^{s-2}}\Big),
$$
which holds for all sufficently large $n$ (without any moment assumption).
Hence, with some constant $c$
\begin{equation}
\int_{T_n}^{4T_n^4} |p_n''(x)\log p_n(x)|\,dx \, \leq \, 
c\int_{T_n}^{4T_n^4} x^2\, |p_n''(x)|\,dx + o\Big(\frac{1}{n^{s-2}}\Big). 
\end{equation}

Now, on the interval $[T_n,4T_n^4]$ we apply Lemma 14.2 with $l=2$ to approximate 
the second derivative. It yields
\bee
\int_{T_n}^{4T_n^4} x^2\, |p_n''(x)|\,dx 
 & \leq &
\int_{T_n}^{+\infty} x^2\, |\varphi_m''(x)|\,dx + 
\int_{T_n}^{+\infty} 
\frac{|\psi_{2,n,1}(x)|}{1+|x|^{s-2}}\,dx \cdot o\Big(\frac{1}{n^{(s-2)/2}}\Big) \\
 & &
 + \ 
\int_{T_n}^{4T_n^4}\frac{1}{1+|x|^{m-2}}\,|\psi_{2,n,2}(x)|\,dx \cdot  \big(O(n^{-(m-1)/2}) + o(n^{-(s-2)})\big).
\ene
Here, as in the proof of Lemma 12.1, the first integral on the right-hand side 
is bounded, up to a constant, by
$$
\int_{T_n}^{+\infty} x^4\varphi(x)\,dx = o\Big(\frac{1}{T_n^{s-3}n^{(s-2)/2}}\Big),
$$
and for the second one, we use Cauchy's inequality to estimate it by
$T_n^{-(s-5/2)}$. Similarly, the last integral is bounded by
$$
2 T_n^2\,\bigg(\int_{-\infty}^{+\infty} \psi_{2,n,2}(x)^2\,dx\bigg)^{1/2}
 \, \leq \, 2T_n^2.
$$
Since $T_n^2$ has a logarithmic growth, we conclude that
$$
\int_{T_n}^{4T_n^4} x^2\, |p_n''(x)|\,dx = 
o\Big(\frac{1}{n^{(s-2)/2}\,(\log n)^{(s-3)/2}}\Big),
$$
so a similar bound also holds for the left integral in (14.5).

To deal with the remaining values of $x$, we will consider the set
$S_1 = \big\{x > 4T_n^4: p_n(x) \leq \frac{1}{2}\,e^{-4\sqrt{x}}\,\big\}$
and its complement $S_2 = (4T_n^4,+\infty) \setminus S_1$.
By Proposition 6.3, for all sufficiently large $n$, and with some 
constants $c,c'$ we have
\bee
\int_{S_1} |p_n''(x)\log p_n(x)|\,dx 
 & \leq & 
c \int_{S_1} \sqrt{p_n(x)}\, |\log p_n(x)|\,dx \\
 & \leq & 
c' \int_{4T_n^4}^{+\infty} \sqrt{x}\, e^{-2\sqrt{x}}\,dx
 \, = \, o\Big(\frac{1}{n^{s-2}}\Big).
\ene
On the other hand, applying (14.2) on the set $S_2$, we get 
\bee
\int_{S_2} |p_n''(x)\log p_n(x)|\,dx|
 & \leq &
c \int_{S_2} |p_n''(x)| \sqrt{x}\,dx \\
 & \leq & 
c' \int_{4T_n^4}^{+\infty} x^{5/2} \varphi(x)\,dx +
c' \int_{4T_n^4}^{+\infty} \frac{dx}{x^{m-1/2}} \cdot 
o\Big(\frac 1{n^{(s-2)/2}}\Big)\\
 & = &
o\Big(\frac{1}{T_n^{2(2m-3)} n^{(s-2)/2}}\Big).
\ene
Combining the two estimates, the theorem is proved.

\vskip5mm
{\bf Remark 14.3.} If $2<s<4$, the expansion (14.1) becomes
\be
I(Z_n||Z) = o\Big(\frac 1{n^{(s-2)/2}\,(\log n)^{(s-3)/2}}\Big).
\en
This formulation does not include the case $s=2$. In case $s>2$, we expect 
that the bound (14.6) may be improved further. However, a possible improvement 
may concern the power of the logarithmic term, only. This can be illustrated 
by means of the example of densities of the form
$$
p(x) = \int_{\sigma_0}^{+\infty} \varphi_\sigma(x)\,dP(\sigma) \qquad (x \in \R),
$$
that is, mixtures of densities of normal distributions on the line with mean
zero, where $P$ is a (mixing) probability measure supported on the half-axis $(\sigma_0,+\infty)$ with $\sigma_0 > 0$. A natural variance constraint
on $P$ is that
\be
\int_{-\infty}^{+\infty} x^2 p(x)\,dx = 
\int_{\sigma_0}^{+\infty} \sigma^2\,dP(\sigma) = 1,
\en
so we should assume that $0 < \sigma_0 < 1$.

First, let us note that, by the convexity of the Fisher information,
$$
I(p) \leq \int_{\sigma_0}^{+\infty} I(\varphi_\sigma)\,dP(\sigma)
= \int_{\sigma_0}^{+\infty} \frac{1}{\sigma^2}\,dP(\sigma) \leq
\frac{1}{\sigma_0^2}, 
$$
hence, $I(p)$ is finite. On the other hand, given $\eta > s/2$,
it is possible to construct the measure $P$ to satisfy (14.7) and with
$$
D(Z_n||Z)\, \geq\, \frac{c}{n^{(s-2)/2} \, (\log n)^\eta},
$$
for all $n$ large enough, and with a constant $c$ depending on $s$ and $\eta$,
only (cf. [B-C-G2]). For example, one may define $P$ 
on the half-axis $[2,+\infty)$ by its density
$$
\frac{dP(\sigma)}{d\sigma} = \frac{c}{\sigma^{s+1} (\log \sigma)^\eta}, 
\qquad \sigma > 2,
$$
and then extend it to any interval $[\sigma_0,2]$ in an arbitrary way so that
to obtain a probability measure satisfying the requirement (14.7). 
Hence, (14.6) is sharp up to a logarithmic factor.

Finally, let us mention that in case $s=2$, $D(Z_n||Z)$ and therefore
$I(Z_n||Z)$ may decay at an arbitrary slow rate.

\end{document}